\documentclass[preprint,11pt,3p]{elsarticle}
\usepackage{amsmath,amsthm,amsfonts,amssymb,arydshln}
\usepackage{xcolor}
\usepackage{algorithm}
\usepackage{multirow}
\usepackage{mathdots}
\usepackage{multicol,graphicx}
\usepackage{mathrsfs}
\usepackage{booktabs}
\usepackage{comment}
\usepackage{subcaption}
\usepackage{threeparttable}
\usepackage{lipsum}
\usepackage{graphicx}
\usepackage{algpseudocode}
\usepackage{tikz-cd}
\algnewcommand{\Input}{\item[\textbf{Input:}]}
\algnewcommand{\Output}{\item[\textbf{Output:}]}

\usepackage{mathrsfs}
\usepackage{sectsty}
\definecolor{astral}{RGB}{46,116,181}
\sectionfont{\color{astral}}
\newtheorem{theorem}{Theorem}[section]
\newtheorem{lemma}[theorem]{Lemma}
\newtheorem{corollary}[theorem]{Corollary}

\newtheorem{definition}[theorem]{Definition}
\newtheorem{example}[theorem]{Example}
\newtheorem{remark}[theorem]{Remark}
\usepackage{hyperref}
\usepackage{scalerel}
\usepackage{stackengine,wasysym}

\usepackage{tikz}
\pgfdeclarelayer{bg} 
\pgfsetlayers{bg,main}

\usepackage{tikz,xcolor}
\definecolor{lime}{HTML}{A6CE39}
\definecolor{lightblue}{rgb}{0.0, 0.0, 0.5}
\DeclareRobustCommand{\orcidicon}{%
	\begin{tikzpicture}
	\draw[lime, fill=lime] (0,0)
	circle [radius=0.16]
	node[white] {{\fontfamily{qag}\selectfont \tiny ID}};
	\draw[white, fill=white] (-0.0625,0.095)
	circle [radius=0.007];
	\end{tikzpicture}
	\hspace{-2mm}
}
\foreach \x in {A, ..., Z}{%
	\expandafter\xdef\csname orcid\x\endcsname{\noexpand\href{https://orcid.org/\csname orcidauthor\x\endcsname}{\noexpand\orcidicon}}
}

\begin{document}
\begin{frontmatter}
\title{
Block Diagonalization and Multiplicativity of Quaternion $Z$-Block Circulant Matrices with Applications
}
\author{Daochang Zhang$^a$, Yue Zhao$^{b}$, Jingqian Li$^c$, Dijana Mosi\'c$^d$}
\address{
$^{a}$ College of Mathematics, Changchun Normal University, Changchun, 130032, P.R. China.
\\\textit{E-mail}: \texttt{daochangzhang@126.com}\\

$^{b}$ Department of Math \& Stat, Georgia State University, Atlanta, GA 30303, USA.
\\\textit{E-mail}: \texttt{yuezhao0303@163.com}\\

$^{c}$ School of Sciences, Tianjin University of Technology, Tianjin, 300384, P.R. China.
\\\textit{E-mail}: \texttt{JingqianLi85@163.com}\\

$^{d}$ Faculty of Sciences and Mathematics, University of Ni\v s, P.O.
Box 224, 18000 Ni\v s, Serbia.
\\\textit{E-mail}: \texttt{dijana@pmf.ni.ac.rs}\\
}

\begin{abstract}
The motivation of this paper is twofold. We firstly investigate the block-diagonalization of the $z$-block circulant matrix $\mathtt{bcirc_z}(\mathcal A)$, and develop the algorithm $\mathtt{bcirc_z}$-inv for computing the inverse of $\mathtt{bcirc_z}(\mathcal A)$.
Second, we establish the multiplicativity of the operator $\mathtt{bcirc_z}$, which yields the equivalence between the QT-product of tensors and the multiplication of their corresponding $z$-block circulant matrices. 

As applications of these results, we conduct large-scale numerical experiments and analyze the scalability of the Tikhonov-regularized model. We also obtain polar decompositions under the QT-product, develop an algorithm for computing them, and present a large-scale test and a scalability analysis. In addition, we evaluate several rotation strategies based on the QT-polar decomposition through video rotation applications. The results show that the decomposition remains stable and maintains inter-frame consistency while accurately preserving colors.
\end{abstract}


\end{frontmatter}

\section{Introduction}
Quaternion algebra, introduced by Hamilton in 1843 to address the problem of representing rotations in a three-dimensional space within mechanics,  generalizes the complex numbers to four dimensions. It forms an associative yet noncommutative algebra. The fundamental multiplication relations,
$\mathbf{i}\mathbf{j}=\mathbf{k}$, $\mathbf{j}\mathbf{k}=\mathbf{i}$, and $\mathbf{k}\mathbf{i}=\mathbf{j}$, naturally encode the algebraic structure of the cross product in $\mathbb{R}^3$. Moreover, quaternions constitute a skew field, preserving the essential algebraic properties that underlie linear and spectral analysis, Fourier transforms, and other foundational constructions in signal processing and applied mathematics.

The expressive power of quaternion algebra has led to extensive applications in signal and image processing~\cite{MFBC2023}, robotics~\cite{CHOU}, computer graphics~\cite{SHOE}, and texture analysis~\cite{KUNZE}. In particular, a quaternion matrix provides an efficient representation of color images, where the three imaginary components correspond to the RGB channels. This compact representation facilitates algebraic operations on color data, motivating the development of quaternion-based matrix and tensor decompositions.

Tensors provide a natural generalization of matrices to higher-order arrays and offer a powerful framework for modeling multidimensional data in signal processing, color image processing, facial recognition, and computer vision~\cite{SIAMMKB2013, SIAMSMS2013, SIAMSL2018, WDPRESS2016, MFBC2023,IETSPLCZ2018}. The theoretical foundations of tensor decomposition were comprehensively reviewed by Kolda and Bader~\cite{SIAMRevKB2009}. With the introduction of new tensor products such as the Einstein product~\cite{AE2007}, the T-product~\cite{LAAKM2011}, and the QT-product \cite{AMLzhang2022}, tensor decompositions defined under these products have attracted increasing research interest.

It is worth noting that most existing studies on tensor decompositions based on various products are formulated over the complex field. Classical tensor decomposition models, including the CANDECOMP/PARAFAC (CP) decomposition~\cite{CRRROL1970}, the Tucker decomposition~\cite{TUCKER1966}, the higher-order singular value decomposition (HOSVD)~\cite{Lathauwer2000}, and the tensor singular value decomposition (t-SVD)~\cite{LAAKM2011} have been extensively investigated. Recently, there have been several new generalizations about the SVD in the quaternion division ring, including the quaternion tensor singular value decomposition (QT-SVD) \cite{AMLzhang2022}, the tensor generalized singular value decomposition (T-QGSVD) \cite{HENG2025} and the constrained quaternion singular value decomposition \cite{wangqinwenR2025}. On the other hand, the QT-QR decomposition was introduced in \cite{LiuAMC2026, NLAAZhangyang2024}. Moreover, the polar decomposition of tensors with the Einstein product was investigated in \cite{Erfani2024}.

Particularly, \cite{LAAKM2011} showed that the normalized Discrete Fourier Transform (DFT) matrix block-diagonalizes any complex block circulant matrix. And \cite{Lund2020} investigated the block circulant matrix preserves multiplication with respect to the T-product over complex field. 

Two key observations can be made.
\begin{enumerate}[(i)]
    \item The above diagonalization result mentioned in \cite{LAAKM2011} is not applicable over the quaternion field. Although both \cite{PanNg2024} and \cite{ZHENGMM} demonstrated that the quaternion block circulant matrix can be block-diagonalized under the T-product, which not keep the similar diagonalization structure. 
    
    \item On the other hand, the multiplicative property established in \cite{Lund2020} remains valid over the quaternion skew field under the T-product, as shown in \cite{HENG2025}.
    
\end{enumerate}
Based on these two observations, this paper aim to derive the corresponding results for quaternion tensors under the QT-product. The specific content is organized as follows. 

In Section~2, we review the relevant definitions and present several key lemmas on third-order quaternion tensors, which form the theoretical foundation for the subsequent sections.

Section~3 introduces the known cross structure of the circulant matrix associated with a quaternion tensor $\mathcal A$. For convenience, we call this matrix the $z$-block circulant matrix of $\mathcal A$, denote it by $\mathtt{bcirc_z}(\mathcal A)$, and explore its block-diagonalization. We then show that the operator $\mathtt{bcirc_z}$ preserves multiplication under the QT-product over the quaternion skew field. This property leads to the computational equivalence between the product of the $z$-block circulant matrices and the QT-product of their corresponding quaternion tensors. 

In Section~4, we propose the QT-polar decomposition, and the corresponding $z$-block circulant matrices decomposition via the results as mentioned in Section 3. Then we present a numerical example to verify the accuracy. 

Section 5 presents a $z$-block circulant structure-accelerated inversion algorithm ($\mathtt{bcirc_z}$-inv) for computing the inverse of $\mathtt{bcirc_z}(\cdot)$. For comparison, the quaternion matrix inverse function \texttt{inv} provided by the \textit{MATLAB} toolbox \textit{QTFM} is employed. Moreover, large-scale tests and scalability analyses of both the inversion algorithm ($\mathtt{bcirc_z}$-inv) and the QT-polar decomposition are provided.

In Section 6, large-scale experiments and scalability analysis of the Tikhonov-regularized model are conducted using Theorem \ref{AB} and the algorithm $\mathtt{bcirc_z}$-inv. Subsequently, a practical example is presented to illustrate video rotation via quaternion tensor polar decomposition, which demonstrates that the QT-polar decomposition
preserves frame-to-frame coherence while maintaining stability and faithful color reproduction.

\section{Preliminaries and key lemmas}
In this section, we present a concise summary of the key definitions and notations used in the subsequent analysis.

For clarity and consistency, throughout this paper, we first introduce the notations that \cite{AMLzhang2022} adopted.

We adopt the conventions that scalars by lowercase letters $a, b, \ldots$, vectors by bold lowercase letters $\mathbf{a}, \mathbf{b}, \ldots$, matrices by uppercase letters $A, B, \ldots$, and third-order tensors by calligraphic letters $\mathcal{A}, \mathcal{B}, \ldots$. In addition, \(A^*\) and \(\mathcal{A}^*\) represent the conjugate transpose of the matrix \(A\) and the third-order tensor \(\mathcal{A}\), separately.
For a third-order tensor \(\mathcal{A}\), we employ MATLAB-style indexing that is \(i\)-th horizontal slice \(\mathcal{A}(i, :, :)\), \(j\)-th lateral slice \(\mathcal{A}(:, j, :)\), and \(k\)-th frontal slice \(\mathcal{A}(:, :, k)\). For concise, we denote the \(k\)-th frontal slice by \(\mathcal A^{(k)} := \mathcal{A}(:, :, k)\).

Let \(\mathbb{R}\) denote the field of real numbers, and \(\mathbb{C}\) the field of complex numbers. The set \(\mathbb{Q}\) of quaternions forms a four-dimensional vector space on \(\mathbb{R}\), with the ordered basis \(\mathbf{1}, \mathbf{i}, \mathbf{j}, \mathbf{k}\). 

A quaternion \(q \in \mathbb{Q}\) can be written in the form
\[
q = q_0 + q_1\mathbf{i} + q_2\mathbf{j} + q_3\mathbf{k},
\]
where \(q_0, q_1, q_2, q_3 \in \mathbb{R}\), and the imaginary units \(\mathbf{i}, \mathbf{j}, \mathbf{k}\) satisfy the following multiplication rules:
\[
\begin{cases}
\mathbf{i}^2 = \mathbf{j}^2 = \mathbf{k}^2 = -1, \\
\mathbf{ij} = -\mathbf{ji} = \mathbf{k}, \\
\mathbf{jk} = -\mathbf{kj} = \mathbf{i}, \\
\mathbf{ki} = -\mathbf{ik} = \mathbf{j}.
\end{cases}
\]
Moreover, the conjugate of a quaternion $q$ is defined as $q^*=q_0 - q_1\mathbf{i} - q_2\mathbf{j} - q_3\mathbf{k}$, and its norm is given by $|q|=\sqrt{q^*q}=\sqrt{q_0^2+ q_1^2+ q_2^2+ q_3^2}$. Note that any quaternion $q$ can be regarded as a two-dimensional vector over complex field $\mathbb C$, and can be uniquely expressed as $$q=c_1+jc_2,$$ where $c_1=q_0 + q_1\mathbf{i}$ and $c_2= q_2- q_3\mathbf{i}$.

Let \(\mathbb{Q}^n\) denote the set of \(n\)-dimensional quaternion vectors. Similarly, let \(Q \in \mathbb{Q}^{m \times n}\) represent the set of \(m \times n\) quaternion matrices. Any quaternion matrix can be expressed in the form
$$Q=Q_{\mathbf{d}}+\mathbf{j}Q_{\mathbf{c}},$$where $Q_{\mathbf{d}}$, $Q_{\mathbf{c}}\in\mathbb C^{m\times n}$. 

Next, we present the definition of eigenvalues for quaternion matrices, following the discussion in~\cite{LAAzfzl1997}.
\begin{definition}
    Given a quaternion matrix \( A \in \mathbb{Q}^{n \times n} \), if there exists a nonzero quaternion vector \( x \in \mathbb{Q}^n \) and a scalar \( \lambda \in \mathbb{Q} \) such that
\[
Ax = x\lambda \quad (\text{respectively, } Ax = \lambda x),
\]
then \( \lambda \) is called a \emph{right} (respectively, \emph{left}) eigenvalue of \( A \), and \( x \) is the corresponding right (left) eigenvector.
\end{definition}
Since quaternion multiplication is non-commutative, eigenvalues of quaternion matrices are classified as left and right eigenvalues, which generally have no direct relationship. In this work, we consider only right eigenvalues.

Furthermore, let \(\mathcal{A} \in \mathbb{Q}^{n_1 \times n_2 \times n_3}\) be a third-order quaternion tensor, represented as 
\begin{align*}
\mathcal{A} =\mathcal A_{\mathbf{d}}+\mathbf{j\mathcal {A}_{\mathbf{c}}},
\end{align*} where $\mathcal A_{\mathbf{d}}$, $\mathcal A_{\mathbf{c}}\in\mathbb C^{n_1\times n_2\times n_3}$.
We define $\mathtt{bcirc}$ as the block circulant matrix associated with \(\mathcal{A}\), of size \(n_1 n_3 \times n_2 n_3\), given by:
\begin{align*}
    \mathtt{bcirc}(\mathcal{A}) :=
\begin{bmatrix}
\mathcal A^{(1)} &\mathcal A^{(n_3)} & \mathcal A^{(n_3-1)} & \cdots & \mathcal A^{(2)} \\
\mathcal A^{(2)} & \mathcal A^{(1)} & \mathcal A^{(n_3)} & \cdots & \mathcal A^{(3)} \\
\vdots & \vdots & \vdots & \ddots & \vdots \\
\mathcal A^{(n_3)} & \mathcal A^{(n_3-1)} & \mathcal A^{(n_3-2)} & \cdots & \mathcal A^{(1)}
\end{bmatrix}
\end{align*}

Moreover, the operation $\mathtt{unfold}$ maps a tensor $\mathcal{A} \in \mathbb{Q}^{n_1 \times n_2 \times n_3}$ into an $n_1 n_3 \times n_2$ block column vector, corresponding to the first block column of $\mathtt{bcirc}(\mathcal A)$. The inverse operation, denoted \(\mathtt{fold}\), reconstructs the original tensor from its unfolded representation. These operations are defined as follows \cite{LAAKM2011}:
\begin{align*}
    \mathtt{unfold}(\mathcal{A}):=\begin{bmatrix}
        \mathcal A^{(1)}\\\mathcal A^{(2)}\\\vdots\\\mathcal A^{(n_3)}
    \end{bmatrix},\quad \mathtt{fold}(\mathtt{unfold}(\mathcal{A})):=\mathcal{A}.
\end{align*}
It is well-known that the normalized Discrete Fourier Transform (DFT) matrix diagonalizes any complex circulant matrix \cite{JHUPgl1996}. Consequently, for any tensor $\mathcal A\in\mathbb{C}^{n_1\times n_2\times n_3 }$, we have 
\begin{align}\label{oldcirc}
    (F_{n_3}\otimes I_{n_1})\mathtt{bcirc}(\mathcal A)(F_{n_3}^*\otimes I_{n_2})=\begin{bmatrix}
        \hat{\mathcal A}^{(1)}&&&\\&\hat{\mathcal A}^{(2)}&&\\&&\ddots&\\&&&\hat{\mathcal A}^{(n_3)}
    \end{bmatrix}:=diag(\hat{\mathcal A}),
\end{align}
where $\hat{\mathcal A}^{(s)}$, $s=1,2,\cdots,n_3$ are frontal slices of $\hat{\mathcal A}$. $F_{n_3}\in\mathbb{C}^{n_3\times n_3}$ is the normalized DFT matrix, where $F_{n_3}(i,j)=n_3^{-1/2}\omega^{(i-1)(j-1)}$  such that  $i,j\in1,2,\cdots, n_3$ and $\omega=e^{\frac{-2\pi \mathbf i}{n_3}}.$ Additionally, for any quaternion tensor $\mathcal A\in \mathbb Q^{n_1\times n_2\times n_3}$, applying the DFT to $\mathcal A$ results in a quaternion tensor $\hat{\mathcal A}\in\mathbb Q^{n_1\times n_2\times n_3},$ where 
\begin{equation}\label{fftq(i,j,:)}
\begin{aligned}
\hat{\mathcal A}(i,j,:)&=\sqrt{n_3}F_{n_3}\mathcal A(i,j,:)=\sqrt{n_3}F_{n_3}\mathcal{A}_{\mathbf{d}}(i,j,:)+\mathbf{j}\sqrt{n_3}P_{n_3}F_{n_3}\mathcal{A}_{\mathbf{c}}(i,j,:)\\&=\hat{\mathcal{A}}_{\mathbf{d}}(i,j,:)+\mathbf{j}P_{n_3}\hat{\mathcal{A}}_{\mathbf{c}}(i,j,:)
\end{aligned}
\end{equation}
such that $i\in 1,2,\cdots,n_1$ and $j\in 1,2,\cdots,n_2$. In addition, the permutation matrix \(P_{n_3} = (P_{ij})\in\mathbb R^{n_3\times n_3}\) is defined as $P_{11}=1$ and
\[
P_{ij} =
\begin{cases}
1, & \text{if } i + j = n_3 + 2, \quad 2 \leq i,~ j \leq n_3, \\
0, & \text{otherwise}.
\end{cases}
\]
The above operation \eqref{fftq(i,j,:)} is represented by a command $fftq(\cdot)$, i.e., $\hat{\mathcal{A}}=fftq(\mathcal{A},[\ ],3)$. Its inverse operation is $\mathcal{A}=ifftq(\hat{\mathcal{A}},[\ ], 3)$.
Then we have the following routine computation: 
\begin{align}\label{zhao2.1}
   \mathtt{unfold}(\hat {\mathcal A})&=\sqrt{n_3}(F_{n_3}\otimes I_{n_1})\mathtt{unfold}(\mathcal A)\text{ and } \mathtt{unfold}(\mathcal A)=\frac{1}{\sqrt{n_3}}(F_{n_3}^*\otimes I_{n_1})\mathtt{unfold}(\hat {\mathcal A})
\end{align}
as also discussed in \cite{AMLzhang2022}. 

In what follows, we present several crucial definitions and lemmas. First, we present the definition of the QT-product as follows:
\begin{definition}\cite{AMLzhang2022}\label{defAB}
Let $\mathcal{A}_{\mathbf{d}},~\mathcal{A}_{\mathbf{c}} \in \mathbb{C}^{n_1 \times r \times n_3}$, 
and consider tensors $\mathcal{A} = \mathcal{A}_{\mathbf{d}} + \mathbf{j}\mathcal{A}_{\mathbf{c}} 
\in \mathbb{Q}^{n_1 \times r \times n_3}$ and 
$\mathcal{B} \in \mathbb{Q}^{r \times n_2 \times n_3}$. 
Here, the symbol~$\otimes$ denotes the Kronecker product. We define
\begin{align*}
    \mathcal A*_Q\mathcal B\stackrel{.}{=}\mathtt{fold}\Big(\big(\mathtt{bcirc}(\mathcal A_{\mathbf{d}})+\mathbf j\mathtt{bcirc}(\mathcal A_{\mathbf{c}})\cdot (P_{n_3}\otimes I_r)\big)\cdot\mathtt{unfold}(\mathcal B)\Big)\in \mathbb{Q}^{n_1 \times n_2 \times n_3}.
    \end{align*}
\end{definition}

We present several basic definitions that will be used frequently in the following sections.

\begin{definition}\cite{AMLzhang2022}
The $n \times n \times n_3$ identity quaternion tensor $\mathcal{I}_{nnn_3}$ is the tensor whose first frontal slice is the identity matrix, and whose remaining frontal slices are zero matrices.
\end{definition}

\begin{definition}\cite{AMLzhang2022}\label{defunitary}
The $n\times n\times n_3$ quaternion tensor $\mathcal{U}$ is said to be unitary if 
\[
\mathcal{U}^**_Q\mathcal{U}
= \mathcal{U}*_Q\mathcal{U}^*
= \mathcal{I}_{n n n_3}.
\]
\end{definition}

\begin{definition}\cite{AMLzhang2022}\label{fdiagonal}
The tensor \(\mathcal{P} \in \mathbb{Q}^{n_1 \times n_2 \times n_3}\) is called an \emph{f-diagonal} quaternion tensor if each of its frontal slices is a diagonal quaternion matrix.
\end{definition}

\begin{definition}\cite{LiuAMC2026}
The tensor \(\mathcal{P} \in \mathbb{Q}^{n \times n \times n_3}\) is called an \emph{f-upper (lower)-triangular} quaternion tensor if each of its frontal slices is an upper (lower)-triangular quaternion matrix.
\end{definition}
Now, the conjugate transpose of a third-order tensor over the complex field is defined as follows.
\begin{definition}\cite{LAAKM2011}
Let $\mathcal{A}\in\mathbb{C}^{n_1\times n_2\times n_3}$. 
Its conjugate transpose $\mathcal{A}^*\in\mathbb{C}^{\,n_2\times n_1\times n_3}$ 
is defined by conjugately transposing each frontal slice of $\mathcal{A}$ 
and reversing the order of slices $2$ through $n_3$.
\end{definition}

Next, we introduce the notion of the conjugate transpose of a third-order quaternion tensor.
\begin{definition}\cite{AMLzhang2022}
The conjugate transpose of a quaternion tensor 
$
\mathcal{A}=\mathcal A_{\mathbf{d}}+\mathbf{j}\mathcal A_{\mathbf{c}}
    \in\mathbb{Q}^{n_1\times n_2\times n_3},
$
denoted by $\mathcal{A}^*\in\mathbb{Q}^{\,n_2\times n_1\times n_3}$, is defined by
\[
\mathtt{unfold}(\mathcal{A}^*)
    = \mathtt{unfold}(\mathcal A_{\mathbf{d}}^*)
      - (P_{n_3}\otimes I_{n_2})\,\mathtt{unfold}(\mathcal A_{\mathbf{c}}^*)\,\mathbf{j}.
\]
\end{definition}

\section{Block-Diagonalization and Multiplicative Equivalence of Quaternion $Z$-Block Circulant Matrices}

In this section, we first introduce the known cross structure of the circulant matrix of the quaternion tensor $\mathcal A$, for convenience, called the $z$-block circulant matrix and denoted by $\mathtt {bcirc_z}(\mathcal A)$, and then investigate its block-diagonalization.
\begin{definition}[\textbf{Z-Block Circulant Matrix}]\label{bcircz}
    Let $\mathcal A\in \mathbb Q^{n_1\times r\times n_3}$, $\mathcal{A}_{\mathbf{d}}\in \mathbb C^{n_1\times r\times n_3}$, and $\mathcal{A}_{\mathbf{c}}\in \mathbb C^{n_1\times r\times n_3}$ satisfy that $$\mathcal A=\mathcal{A}_{\mathbf{d}} + \mathbf{j} \mathcal{A}_{\mathbf{c}}\in \mathbb {Q}^{n_1\times r\times n_3}.$$ 
    Then the z-block circulant matrix of $\mathcal A$ is defined by $\mathtt{bcirc_z}(\mathcal A)$ as follows:
    \begin{align}\label{defbcirc_z}
    \mathtt{bcirc_z}(\mathcal A)=\mathtt{bcirc}(\mathcal{A}_{\mathbf{d}})+\mathbf{j}\mathtt{bcirc}(\mathcal{A}_{\mathbf{c}})(P_{n_3}\otimes I_r),
    \end{align}
    where $P_{n_3}$ is as defined in Definition \ref{defAB}.
\end{definition}

Based on Definition \ref{bcircz}, we conclude that the explicit structure of the $z$-block circulant matrix $\mathtt{bcirc_z}(\mathcal A)$ is as follows:
    \begin{align*}
        \mathtt{bcirc_z}(\mathcal A)=
        \begin{bmatrix}
        \mathcal{A}_{\mathbf{d}}^{(1)} + \mathbf{j}\mathcal{A}_{\mathbf{c}}^{(1)} & \mathcal{A}_{\mathbf{d}}^{(n_3)} +  \mathbf{j}\mathcal{A}_{\mathbf{c}}^{(2)} & \cdots & \mathcal{A}_{\mathbf{d}}^{(2)} + \mathbf{j}\mathcal{A}_{\mathbf{c}}^{(n_3)} \\
        \mathcal{A}_{\mathbf{d}}^{(2)} + \mathbf{j}\mathcal{A}_{\mathbf{c}}^{(2)} & \mathcal{A}_{\mathbf{d}}^{(1)} + \mathbf{j}\mathcal{A}_{\mathbf{c}}^{(3)} & \cdots & \mathcal{A}_{\mathbf{d}}^{(3)} + \mathbf{j}\mathcal{A}_{\mathbf{c}}^{(1)} \\
        \mathcal{A}_{\mathbf{d}}^{(3)} + \mathbf{j}\mathcal{A}_{\mathbf{c}}^{(3)} & \mathcal{A}_{\mathbf{d}}^{(2)} + \mathbf{j}\mathcal{A}_{\mathbf{c}}^{(4)} & \cdots & \mathcal{A}_{\mathbf{d}}^{(4)} + \mathbf{j}\mathcal{A}_{\mathbf{c}}^{(2)} \\
        \vdots & \vdots & \ddots & \vdots \\
        \mathcal{A}_{\mathbf{d}}^{(n_3)} + \mathbf{j}\mathcal{A}_{\mathbf{c}}^{(n_3)} & \mathcal{A}_{\mathbf{d}}^{(n_3-1)} + \mathbf{j}\mathcal{A}_{\mathbf{c}}^{(1)} & \cdots & \mathcal{A}_{\mathbf{d}}^{(1)} + \mathbf{j}\mathcal{A}_{\mathbf{c}}^{(n_3-1)}
        \end{bmatrix}.
    \end{align*}


    
In what follows, we provide a concise and transparent relation between $\mathtt{bcirc_z}(\mathcal{A})$ and $\operatorname{diag}(\hat{\mathcal{A}})$ in the quaternion division ring as follows. 

\begin{theorem}\label{zdc}
Let $\mathcal{A}_{\mathbf{d}}$ and $\mathcal{A}_{\mathbf{c}}\in \mathbb C^{n_1\times n_2\times n_3}$ satisfy that $\mathcal A=\mathcal{A}_{\mathbf{d}} + \mathbf{j} \mathcal{A}_{\mathbf{c}}\in\mathbb{Q}^{{n_1\times n_2\times n_3}}.$ Denote that $\hat{\mathcal A}$ is the DFT of $\mathcal{A}$. Then, 
\begin{align}\label{bcircdiag}
    (F_{n_3}\otimes I_{n_1})\mathtt{bcirc_z}(\mathcal A)(F_{n_3}^*\otimes I_{n_2})=diag(\hat{\mathcal A}).
\end{align}
\begin{proof}
        Recall that the DFT matrix $F_{n_3}$ satisfies $F_{n_3}F_{n_3}=P_{n_3},$ where $P_{n_3}$ is a permutation matrix. Moreover, the conjugate of $F_{n_3}$, denoted $\bar{F}_{n_3}$, satisfies $\bar{F}_{n_3}=F_{n_3}^{*}.$ Then, we conclude that $$F_{n_3}\mathbf{j}=\mathbf{j}\bar{F}_{n_3}=\mathbf{j}F_{n_3}^{*}=\mathbf{j}F_{n_3}^{*}(F_{n_3}F_{n_3}^*)=\mathbf{j}P_{n_3}F_{n_3}.$$
        In addition, the following equations are helpful for computation:$$F_{n_3}P_{n_3}=P_{n_3}F_{n_3}=F_{n_3}^*,\quad F_{n_3}^*P_{n_3}=P_{n_3}F_{n_3}^*=F_{n_3}.$$
Substituting \eqref{defbcirc_z} into the left-hand side of \eqref{bcircdiag}, we obtain
\begin{align}\label{equ1}
\notag&(F_{n_3}\otimes I_{n_1})\Big(\mathtt{bcirc_z}(\mathcal A)\Big)(F_{n_3}^*\otimes I_{n_2})
            \\\notag=&(F_{n_3}\otimes I_{n_1})\Big(\mathtt{bcirc}(\mathcal A_{\mathbf{d}})+\mathbf{j}\mathtt{bcirc}(\mathcal A_{\mathbf{c}})(P_{n_3}\otimes I_{n_2})\Big)(F_{n_3}^*\otimes I_{n_2})
            \\=&(F_{n_3}\otimes I_{n_1})\mathtt{bcirc}(\mathcal A_{\mathbf d})(F_{n_3}^*\otimes I_{n_2})
            +(F_{n_3}\otimes I_{n_1})\mathbf{j}\mathtt{bcirc}(\mathcal A_{\mathbf{c}})(F_{n_3}\otimes I_{n_2})
\end{align}
It is worth noting that the matrix-form identity $(UV)\otimes(YW)=(U\otimes Y)(V\otimes W)$. Moreover, we observe that $(F_{n_3}\otimes I_{n_1})\mathbf{j}=(F_{n_3}\mathbf{j}\otimes I_{n_1})=(\mathbf{j}F_{n_3}^*\otimes I_{n_1})=\mathbf{j}(F_{n_3}^*\otimes I_{n_1})=\mathbf{j}(F_{n_3}P_{n_3}\otimes I_{n_1})$, then \eqref{equ1} can be reformulate as 
\begin{align*}
    \notag&(F_{n_3}\otimes I_{n_1})\Big(\mathtt{bcirc_z}(\mathcal A)\Big)(F_{n_3}^*\otimes I_{n_2})\\&=(F_{n_3}\otimes I_{n_1})\mathtt{bcirc}(\mathcal A_{\mathbf{d}})(F_{n_3}^*\otimes I_{n_2})
            +\mathbf{j}(F_{n_3}\otimes I_{n_1})(P_{n_3}\otimes I_{n_1})\mathtt{bcirc}(\mathcal A_{\mathbf{c}})(F_{n_3}^*\otimes I_{n_2})(P_{n_3}\otimes I_{n_2}).
\end{align*}
Therefore, combining \eqref{oldcirc} and \eqref{fftq(i,j,:)}, the above equation can be rewritten as
\begin{align*}
    (F_{n_3}\otimes I_{n_1})\Big(\mathtt{bcirc_z}(\mathcal A)\Big)(F_{n_3}^*\otimes I_{n_2})=&diag(\hat{\mathcal A}_{\mathbf{d}})+\mathbf j(P_{n_3}\otimes I_{n_1})diag(\hat{\mathcal A_{\mathbf{c}}})(P_{n_3}\otimes I_{n_2})=diag(\hat{\mathcal A}).
\end{align*}        
Hence, the proof is completed.
\end{proof}
\end{theorem}

We now proceed to define the following two inverse operations.

The operator $diag^{-1}(\cdot)$ is denoted by:
$$diag^{-1}(\begin{bmatrix}
        \hat{\mathcal A}_{(1)}&&&\\&\hat{\mathcal A}_{(2)}&&\\&&\ddots&\\&&&\hat{\mathcal A}_{(n_3)}
    \end{bmatrix})=\hat{\mathcal A}.$$

The operator $\mathtt{bcirc_z}^{-1}(\cdot)$ is defined as:
  \begin{align}\label{inverse operator2}
        \mathcal A=\mathtt{bcirc_z}^{-1}(
        \begin{bmatrix}
        \mathcal{A}_{\mathbf{d}}^{(1)} + \mathbf{j}\mathcal{A}_{\mathbf{c}}^{(1)} & \mathcal{A}_{\mathbf{d}}^{(n_3)} +  \mathbf{j}\mathcal{A}_{\mathbf{c}}^{(2)} & \cdots & \mathcal{A}_{\mathbf{d}}^{(2)} + \mathbf{j}\mathcal{A}_{\mathbf{c}}^{(n_3)} \\
        \mathcal{A}_{\mathbf{d}}^{(2)} + \mathbf{j}\mathcal{A}_{\mathbf{c}}^{(2)} & \mathcal{A}_{\mathbf{d}}^{(1)} + \mathbf{j}\mathcal{A}_{\mathbf{c}}^{(3)} & \cdots & \mathcal{A}_{\mathbf{d}}^{(3)} + \mathbf{j}\mathcal{A}_{\mathbf{c}}^{(1)} \\
        \mathcal{A}_{\mathbf{d}}^{(3)} + \mathbf{j}\mathcal{A}_{\mathbf{c}}^{(3)} & \mathcal{A}_{\mathbf{d}}^{(2)} + \mathbf{j}\mathcal{A}_{\mathbf{c}}^{(4)} & \cdots & \mathcal{A}_{\mathbf{d}}^{(4)} + \mathbf{j}\mathcal{A}_{\mathbf{c}}^{(2)} \\
        \vdots & \vdots & \ddots & \vdots \\
        \mathcal{A}_{\mathbf{d}}^{(n_3)} + \mathbf{j}\mathcal{A}_{\mathbf{c}}^{(n_3)} & \mathcal{A}_{\mathbf{d}}^{(n_3-1)} + \mathbf{j}\mathcal{A}_{\mathbf{c}}^{(1)} & \cdots & \mathcal{A}_{\mathbf{d}}^{(1)} + \mathbf{j}\mathcal{A}_{\mathbf{c}}^{(n_3-1)}
        \end{bmatrix}).
    \end{align}

\begin{remark}
    Let
$
\mathcal{A} = \mathcal{A}_{\mathbf{d}} + \mathbf{j}\mathcal{A}_{\mathbf{c}} 
    \in \mathbb{Q}^{n_1 \times n_2 \times n_3},
$
where $\mathcal{A}_{\mathbf{d}}, \mathcal{A}_{\mathbf{c}} \in \mathbb{C}^{n_1 \times n_2 \times n_3}$. 
The tensor $\hat{\mathcal{A}} \in \mathbb{Q}^{n_1 \times n_2 \times n_3}$ is obtained by applying the DFT to 
$\mathcal{A}_{\mathbf{d}}$ and $\mathcal{A}_{\mathbf{c}}$ along the third mode. Then
$$fftq(\mathcal{A},[\ ],3)=\hat{\mathcal A}=idiag((F_{n_3}\otimes I_{n_1})\mathtt{bcirc_z}(\mathcal A)(F_{n_3}^*\otimes I_{n_2})).$$
$$ ifftq(\hat{\mathcal A},[\ ],3)=\mathcal A=\mathtt{ibcirc_z}((F_{n_3}^*\otimes I_{n_2})diag(\hat{\mathcal A})(F_{n_3}\otimes I_{n_1})).$$
\end{remark}

Next, we present our second main result, which establishes that the operator
$\mathtt{bcric_z}$ preserves multiplication with respect to the QT-product. This result shows that $\mathtt{bcric_z}$ defines an algebra homomorphism from the QT-product algebra of quaternion tensors to the algebra of structured matrices under matrix multiplication, thereby transforming noncommutative quaternion tensor operations into matrix operations.
\begin{theorem}\label{AB}
    Let $\mathcal A\in \mathbb Q^{n_1\times r\times n_3}$ and $\mathcal B\in \mathbb Q^{r\times n_2\times n_3}.$  Then $\mathtt{bcric_z}(\mathcal A*_Q\mathcal B)=\mathtt{bcric_z}(\mathcal A)\mathtt{bcric_z}(\mathcal B)$.
\end{theorem}
\begin{proof}
        From the definitions of QT-product and $\mathtt{bcirc_z}$, we have
        \begin{align*}
            \mathcal A*_Q\mathcal B=\mathtt{fold}(\mathtt{bcirc_z(\mathcal A)\mathtt{unfold}(\mathcal B)}).
        \end{align*}
        Furthermore, taking $\mathtt{unfold}$ operation on both sides, we obatin
        \begin{align}\label{dc1}
            \mathtt{unfold}(\mathcal A*_Q\mathcal B)\notag&=\mathtt{bcirc_z(\mathcal A)\mathtt{unfold}(\mathcal B)}
            \\&=\frac{1}{\sqrt{n_3}}\mathtt{bcirc_z}(\mathcal A)(F_{n_3}^*\otimes I_r)\mathtt{unfold}(\hat{\mathcal B}),
        \end{align}
        and the second equality is derived from \eqref{zhao2.1}. In addition, we can still get the following equality by utilizing \eqref{zhao2.1},
        \begin{align*}
            \mathtt{unfold}(\widehat{\mathcal A*_Q\mathcal B})&=\sqrt{n_3}(F_{n_3}\otimes I_{n_1})\mathtt{unfold}(\mathcal A*_Q\mathcal B)
            \\&=(F_{n_3}\otimes I_{n_1})\mathtt{bcirc_z(\mathcal A)}(F_{n_3}^*\otimes I_r)\mathtt{unfold}(\hat{\mathcal B}),
        \end{align*}
        where the second equality is obtained by substituting $\mathtt{unfold}(\mathcal A*_Q\mathcal B)$ with \eqref{dc1}. 
        By utilizing Theorem \ref{zdc}, it directly follows that\begin{align*}
            &(F_{n_3}\otimes I_{n_1})\mathtt{bcirc_z}(\mathcal A*_Q\mathcal B)(F_{n_3}^*\otimes I_{n_2})\\&=(F_{n_3}\otimes I_{n_1})\mathtt{bcirc_z}(\mathcal A)(F_{n_3}^*\otimes I_{r})(F_{n_3}\otimes I_{r})\mathtt{bcirc_z}(\mathcal B)(F_{n_3}^*\otimes I_{n_2}),
        \end{align*}
        which implies \begin{align*}
            \mathtt{bcirc_z}(\mathcal A*_Q\mathcal B)=\mathtt{bcirc_z}(\mathcal A)\mathtt{bcirc_z}(\mathcal B).
        \end{align*}
        Hence, the proof is complete.
    \end{proof}
Next, we establish the multiplicative equivalence between $z$-block circulant matrices and their corresponding quaternion tensors under the QT-product as follows. Applying Theorem \ref{AB} and the inverse operator $\mathtt{bcirc_z}^{-1}(\cdot)$ in \eqref{inverse operator2}, we obtain the following corollary.
\begin{theorem}\label{newrelation}
    Let $\mathcal A\in \mathbb Q^{n_1\times r\times n_3},$ $\mathcal B\in \mathbb Q^{r\times n_2\times n_3},$ and $\mathcal C\in \mathbb Q^{n_1\times n_2\times n_3}.$  Then
    \begin{align*}
       \mathcal A*_Q\mathcal B=\mathcal C\Longleftrightarrow \mathtt{bcric_z}(\mathcal A)\mathtt{bcric_z}(\mathcal B)=\mathtt{bcric_z}(\mathcal C).
    \end{align*}
    \end{theorem}
    
    \begin{remark}
    It is worth mentioning that \cite{AMLzhang2022} discussed the following relationship when $\mathcal A\in\mathbb{Q}^{n_1\times r\times n_3},\mathcal B\in\mathbb{Q}^{r\times n_2\times n_3},$ and $\mathcal C\in\mathbb{Q}^{n_1\times n_2\times n_3}$ under the QT-product:
\begin{align*}
    \mathcal A*_Q\mathcal B=\mathcal C\Longleftrightarrow diag(\hat{\mathcal A})diag(\hat{\mathcal B})=diag(\hat{\mathcal C}),
    \end{align*} where $\hat{\mathcal{A}}$, $\hat{\mathcal B}$, and $\hat{\mathcal{C}}$ are the DFT of 
    $\mathcal{A}$, $\mathcal B$, and $\mathcal{C},$ respectively.
    
    Based on Theorem \ref{zdc}, it follows immediately that 
\[
\begin{tikzcd}[column sep=1em, row sep=1em]
\mathtt{bcirc_z}(\mathcal C)
=
\mathtt{bcirc_z}(\mathcal A)\,
\mathtt{bcirc_z}(\mathcal B)
\arrow[rr, Leftrightarrow]
&&
diag(\hat{\mathcal C})
=
diag(\hat{\mathcal A})
diag(\hat{\mathcal B})
\arrow[dl, Leftrightarrow]
\\
&
\mathcal C
=
\mathcal A *_Q \mathcal B
\arrow[ul, Leftrightarrow]
&
\end{tikzcd}
\]
However, we can not obtain the above conclusion without Theorem \ref{zdc}.
    \end{remark}

    
     

    In addition, we discuss the conjugate transpose equivalence between the quaternion tensor $\mathcal A$ and the $z$-block circulant matrix $\mathtt{bcirc_z}(\mathcal A)$.
    \begin{theorem}\label{ct}
    Let $\mathcal A\in\mathbb Q^{n\times n\times n_3}$ and the conjugate transpose of $\mathcal A$ is defined by $\mathcal A^*$. Then 
    \begin{align*}
        \mathtt{bcirc_z}(\mathcal A^*)=\mathtt{bcirc_z}(\mathcal A)^*.
    \end{align*}
    \begin{proof}
        Observe that the structure of $\mathcal A^*$ is as follows:
        \[
        \mathtt{unfold}({\mathcal A^*})=\mathtt{unfold}(\mathcal A_{\mathtt{d}}^*)-(P_{n_3}\otimes I_n)\mathtt{unfold}(\mathcal A_{\mathbf{c}}^*)\mathbf{j}.
        \]
        Then, by Definition \ref{bcircz}, we provide the $z$-block circulant matrix of $\mathcal A^*$
        \begin{align}\label{Act}
        \mathtt{bcirc_z}(\mathcal A^*)=\mathtt{bcirc}(\mathcal{A}_{\mathbf{d}}^*)+\mathbf{j}\mathtt{bcirc}(-\mathcal{A}_{\mathbf{c}}^{T})(P_{n_3}\otimes I_n),
        \end{align}
        where $\mathcal A_{\mathbf{c}}^T$ is the transpose of $\mathcal A_{\mathbf{c}}$. Moreover, we can rewrite \eqref{Act} as follows:
        \begin{align}\label{reAct}
            \mathtt{bcirc_z}(\mathcal A^*)&\notag=\mathtt{bcirc}(\mathcal{A}_{\mathbf{d}}^*)-\mathtt{bcirc}(\mathbf{j}\mathcal{A}_{\mathbf{c}}^{T})(P_{n_3}\otimes I_n)
            \\\notag&=\mathtt{bcirc}(\mathcal{A}_{\mathbf{d}}^*)-\mathtt{bcirc}\Big(\mathtt{fold}\big((P_{n_3}\otimes I_n)\mathtt{unfold}(\mathcal A_{\mathbf{c}}^*)\big)\Big)\mathbf{j}(P_{n_3}\otimes I_n)
            \\&=\mathtt{bcirc}(\mathcal{A}_{\mathbf{d}}^*)-(P_{n_3}\otimes I_n)\mathtt{bcirc}(\mathcal A_{\mathbf{c}}^*)\mathbf{j},
        \end{align} where the third equality is obtained from 
        \[
        \mathtt{bcirc}\Big(\mathtt{fold}\big((P_{n_3}\otimes I_n)\mathtt{unfold}(\mathcal A_{\mathbf{c}}^*)\big)\Big)=(P_{n_3}\otimes I_n)\mathtt{bcirc}(\mathcal A_{\mathbf{c}}^*)(P_{n_3}\otimes I_n).
        \]
        Additionally, from Definition \ref{bcircz}, we conclude that
        \begin{align}\label{bAct}
            \mathtt{bcirc_z}(\mathcal A)^*\notag&=\mathtt{bcirc}(\mathcal{A}_{\mathbf{d}})^*-(P_{n_3}\otimes I_n)\mathtt{bcirc}(\mathcal A_{\mathbf{c}})^*\mathbf{j}
            \\&=\mathtt{bcirc}(\mathcal{A}_{\mathbf{d}}^*)-(P_{n_3}\otimes I_n)\mathtt{bcirc}(\mathcal A_{\mathbf{c}}^*)\mathbf{j},
        \end{align}
        and the second equality follows from the structure of the block circulant matrix, that is $\mathtt{bcirc}(\mathcal A_{\mathbf{d}}^*)=\mathtt{bcirc}(\mathcal A_{\mathbf{d}})^*$. Therefore, it is directly verified that \eqref{reAct} coincides with \eqref{bAct}, thereby completing the proof.
    \end{proof}
    \end{theorem}

    \begin{remark}
    We state the following identity:
        \begin{align*}
        \mathtt{bcirc_z}(\mathcal A^*)=\mathtt{bcirc_z}(\mathcal A)^*
        \Longleftrightarrow diag(\widehat{\mathcal A^*})=diag(\hat{\mathcal A})^*.
    \end{align*}
    However, $diag(\widehat{\mathcal A^*})\not =diag(\hat{\mathcal A}^*).$
    \end{remark}

    Then we discuss the unitary relation between the quaternion tensor $\mathcal A$ and the $z$-block circulant matrix $\mathtt{bcirc_z}(\mathcal A)$.
    \begin{corollary}\label{unitary}
        A quaternion tensor $\mathcal A\in\mathbb{Q}^{n\times n\times n_3}$ is unitary if and only if $\mathtt{bcirc_z}(\mathcal A)$ is unitary.
        \begin{proof}
            Suppose that $\mathtt{bcirc_z}(\mathcal A)$ is unitary. Then it follows that 
            \[\mathtt{bcirc_z}(\mathcal A)\mathtt{bcirc_z}(\mathcal A)^*=\mathtt{bcirc_z}(\mathcal A)^*\mathtt{bcirc_z}(\mathcal A)=I_{nn_3}.\] By utilizing Theorem \ref{ct} and Theorem \ref{newrelation}, we obtain
            \begin{align*}
            &\mathtt{bcirc_z}(\mathcal A)\mathtt{bcirc_z}(\mathcal A^*)=\mathtt{bcirc_z}(\mathcal A^*)\mathtt{bcirc_z}(\mathcal A)=I_{nn_3}
            \\\Longleftrightarrow\\&\mathcal A*_Q\mathcal A^*=\mathcal A^**_Q\mathcal A=\mathcal I_{nnn_3}.
            \end{align*}
            Thus, the condition is necessary and sufficient.
        \end{proof}
    \end{corollary}

    \begin{remark}
        It is worth noting that, based on Corollary \ref{unitary}, it follows immediately that 
        \begin{align*}
             &\mathtt{bcirc_z}(\mathcal A)\mathtt{bcirc_z}(\mathcal A)^*=\mathtt{bcirc_z}(\mathcal A)^*\mathtt{bcirc_z}(\mathcal A)=I_{nn_3}&\\
             \Longleftrightarrow&\mathcal {A}*_Q\mathcal A^*=\mathcal A^**_Q\mathcal A=\mathcal I_{nnn_3} \\\Longleftrightarrow &diag(\hat{\mathcal A})diag(\hat{\mathcal A})^*=diag(\hat{\mathcal A})^*diag(\hat{\mathcal A})=I_{nn_3}.
        \end{align*}
        where the second necessary and sufficient condition was investigated in \cite{AMLzhang2022} using a different approach. However, we can not obtain the above conclusion without Theorem \ref{ct}.
    \end{remark}

It is well known that a quaternion tensor $\mathcal{A}\in\mathbb{Q}^{n\times n\times n_3}$ is Hermitian if and only if $\mathcal A^*=\mathcal A$, or equivalently, $\mathtt{unfold}(\mathcal A)=\mathtt{unfold}(\mathcal A^{*}).$ Then we can provide the following theorem directly.

\begin{theorem}\label{Hermitian}
Let $\mathcal H\in\mathbb{Q}^{n\times n\times n_3}$ and let $\mathtt{bcirc_z}(\mathcal H)$ denote the $z$-block circulant matrix of $\mathcal H$. Then
\begin{align*}
   \mathcal H \text{ is Hermitian}\Longleftrightarrow \mathtt{bcirc_z({\mathcal H})}
   \text{ is Hermitian}.
\end{align*}

\begin{proof}
We first show the sufficiency.

By Theorem \ref{ct}, if 
$
\mathtt{bcirc_z}(\mathcal{H})
$ is Hermitian,
then
\[\mathtt{bcirc_z}(\mathcal{H}^*)=
\mathtt{bcirc_z}(\mathcal{H})^*
    =\mathtt{bcirc_z}(\mathcal{H}).
\]
It is obvious that $\mathtt{bcirc_z}$ is an invertible operator, i.e.,
\[
\mathcal{H}^* =\mathtt{bcirc_z}^{-1}\big(\mathtt{bcirc_z}(\mathcal{H})^{*}\big) =\mathtt{bcirc_z}^{-1}\big(\mathtt{bcirc_z}(\mathcal{H})\big)
=\mathcal{H}.
\]

We now prove the necessity.

If $\mathcal{H}$ is Hermitian, by Theorem \ref{ct}, we obtain
\[
\mathtt{bcirc_z}(\mathcal{H})^*=\mathtt{bcirc_z}(\mathcal{H}^*)
    =\mathtt{bcirc_z}(\mathcal{H}).
\]
Combining the two directions, the statement follows. 
\end{proof}
\end{theorem}

Moreover, we need to define the inverse of a quaternion tensor under the QT-product and then investigate its equivalence with the inverse of the corresponding $z$-block circulant matrix.

\begin{definition}
Let $\mathcal{P}\in\mathbb{Q}^{n\times n\times n_3}$. The quaternion tensor $\mathcal{U}$ is said to be invertible with respect to the QT-product if there exists a quaternion tensor $\mathcal{Q}\in\mathbb{Q}^{n\times n\times n_3}$ such that
\[
\mathcal{P} *_Q \mathcal{Q}
=
\mathcal{Q} *_Q \mathcal{P}
=
\mathcal{I}_{n n n_3},
\]
where $\mathcal{I}_{n n n_3}$ denotes the identity quaternion tensor under the QT-product. In this case, $\mathcal{Q}$ is called the inverse of~$\mathcal{P}$ and is denoted by $\mathcal{P}^{-1}$. 
\end{definition}

 \begin{lemma}\label{exchangeinverse}
    Let $\mathcal P\in\mathbb{Q}^{n\times n\times n_3}$. Then $\mathcal P$ is invertible if and only if $\mathtt{bcirc_z}(\mathcal P)$ is invertible.
    
    \begin{proof}
        If $\mathcal P$ is invertible, it follows that $\mathcal P*\mathcal P^{-1}=\mathcal I$, by Theorem \ref{newrelation}, we obtain
    \begin{align*}
     \mathtt{bcirc_z}(\mathcal P)\mathtt{bcirc_z}(\mathcal P^{-1})=\mathtt{bcirc_z}(\mathcal I)=I_{nn_3},
     \end{align*}
     which implies $\mathtt{bcirc_z}(\mathcal P)$ is invertible and
     $$\mathtt{bcirc_z}(\mathcal P)^{-1}=\mathtt{bcirc_z}(\mathcal P^{-1}).$$
     Suppose that $\mathtt{bcirc_z}(\mathcal P)$ is invertible, then there exist a quaternion matrix $\mathtt{bcirc_z}(\mathcal Q)$ such that
     \begin{align*}
         \mathtt{bcirc_z}(\mathcal P)\mathtt{bcirc_z}(\mathcal Q)=I_{nn_3}=\mathtt{bcirc_z}(\mathcal I).
     \end{align*}
     It directly follows that $\mathcal P*_Q\mathcal Q=\mathcal I$ by utilizing Theorem \ref{newrelation}. Hence, $\mathcal P$ is said to be an invertible quaternion tensor and denoted by
     $\mathcal P^{-1}=\mathcal Q$.
    \end{proof}
    \end{lemma}

\section{QT-(Right) Polar decomposition}

This section extends the polar decomposition of quaternion matrices to quaternion tensors under the QT-product. We first introduce the notions of \(f\)-positive definite and \(f\)-positive semidefinite Hermitian quaternion tensors, motivated by \cite[Corollaries 6.1 and 6.2]{LAAzfzl1997}. Based on these definitions, we then establish the polar decomposition of quaternion tensors with respect to the QT-product.

\begin{definition}\cite{LAAzfzl1997}\label{MATRIXHP}
An $n\times n$ Hermitian matrix $A$ is said to be \emph{positive (semidefinite)} 
if $x^*Ax >(\ge) 0$ for all nonzero column vectors $x$ of $n$ quaternion components.
\end{definition}

By Definition \ref{MATRIXHP}, we present the following definition. 
\begin{definition}[\textbf{F-Positive Definite (Semidefinite)}]\label{psdef}
    A Hermitian quaternion tensor $\mathcal A\in\mathbb{Q}^{n\times n\times n_3}$ is called f-positive definite (semidefinite) if each frontal slice of $\mathcal A$ is positive definite (semidefinite).
\end{definition}

\begin{remark}\label{NHPSD}\cite{LAAzfzl1997}
Let $A\in\mathbb{Q}^{n\times n}$ with standard eigenvalues $h_1 + k_1\mathbf{i},\,\ldots,\,h_n + k_n\mathbf{i}$ such that the right spectrum
$\sigma_r(A) = [h_1 + k_1\mathbf{i}] \,\cup\, \cdots \,\cup\, [h_n + k_n\mathbf{i}].$
$A$ is normal if and only if there exists a unitary matrix $U$ such that
\[
U^* A U = \operatorname{diag}\{h_1 + k_1\mathbf{i},\, h_2 + k_2\mathbf{i},\, \ldots,\, h_n + k_n\mathbf{i}\}=D,
\]
and $A$ is Hermitian if and only if $k_i=0$ for $1\leq i \leq n$.

Moreover, $A$ is Hermitian positive definite (semidefinite) if and only if  $h_i>0$ ($h_i\geq0$) and $k_i=0$ for $1\leq i \leq n$. 
\end{remark}

It is worth noting that the left polar decomposition for quaternion matrices was presented in \cite[Theorem 7.1]{LAAzfzl1997}. Utilizing the argument used in
\cite[Theorem~7.1]{LAAzfzl1997}, one can analogously show that there exist a
unitary quaternion matrix~$U$ and a Hermitian positive semidefinite
quaternion matrix~$H$ such that the right polar decomposition of a
quaternion matrix~$A$ can be expressed as $A=UH.$
In what follows, we present the next theorem.
\begin{theorem}[\textbf{QT-(Right) Polar decomposition}]\label{QTPOLAR}
    Let $\mathcal{A}\in\mathbb{Q}^{n\times n\times n_3}$. Then there exist a unitary quaternion tensor $\mathcal U$, and a Hermitian and f-positive semidefinite quaternion tensor $\mathcal H$ such that
    \begin{align*}
        \mathcal A=\mathcal U*_Q\mathcal H.
    \end{align*} 
In the meantime,  
\begin{align*}
            \mathtt{bcirc_z}(\mathcal A)=\mathtt{bcirc_z}(\mathcal U)\mathtt{bcirc_z}(\mathcal H),
        \end{align*} where $\mathtt{bcirc_z}(\mathcal U)$ is unitary, and $\mathtt{bcirc_z}(\mathcal H)$ is Hermitian and positive semidefinite. 
    
    \begin{proof}
        From the right polar decomposition of a quaternion matrix, we obtain $\hat{\mathcal A}^{(i)}=\hat{\mathcal U}^{(i)}\hat{\mathcal H}^{(i)}$, where $\hat{\mathcal H}^{(i)}$ is Hermitian and positive semidefinite, and $\hat{\mathcal U}^{(i)}$ is unitary. Moreover, for each $1\le i\le n_{3}$ and every nonzero vector $x^{(i)}\in\mathbb{Q}^{n}$, we have
\[
(x^{(i)})^{*}\hat{\mathcal H}^{(i)}x^{(i)} \ge 0.
\] 
        Notice that $\hat{\mathcal A}^{(i)}$, $\hat{\mathcal H}^{(i)}$, and $\hat{\mathcal U}^{(i)}$ denote the $i$-th frontal slice of $\hat{\mathcal A}$, $\hat{\mathcal H}$, and $\hat{\mathcal U}$, respectively. Then we have
        \begin{align}\label{diagpolar}
            diag(\hat{\mathcal A})=diag(\hat{\mathcal U})diag(\hat{\mathcal H}),
        \end{align} where 
        $diag(\hat{\mathcal U})$ is unitary and $diag(\hat{\mathcal H})$ is Hermitian with 
        \begin{align}\label{diagpsd}
        x^*diag(\hat{\mathcal H})x\geq0,\quad x=
        \begin{bmatrix}
      x^{(1)}&x^{(2)}&\dots&x^{(n_3)}
        \end{bmatrix}^{T},
        \end{align} which implies $diag(\hat{\mathcal H})$ is positive semidefinite.

By Theorem \ref{zdc}, we have
         \begin{align*}
           \mathtt{bcirc_z}(\mathcal U)=(F_{n_3}^*\otimes I_n)diag(\hat{\mathcal U})(F_{n_3}\otimes I_n), 
        \end{align*} then $\mathtt{bcirc_z}(\mathcal U)$ is unitary. Furthermore, by Corollary \ref{unitary}, we obtain $\mathcal U$ is unitary.

        By applying Theorem \ref{zdc} and equation \eqref{diagpolar}, we derive
        \begin{align*}\label{bcirczpolar}
            \mathtt{bcirc_z}(\mathcal A)=\mathtt{bcirc_z}(\mathcal U)\mathtt{bcirc_z}(\mathcal H).
        \end{align*}
         Furthermore, utilizing Theorem \ref{zdc}, we observe that 
        \begin{align*}
           \mathtt{bcirc_z}(\mathcal H)=(F_{n_3}^*\otimes I_n)diag(\hat{\mathcal H})(F_{n_3}\otimes I_n), 
        \end{align*}
        it follows that
        \begin{align*}
            \mathtt{bcirc_z}(\mathcal H)^*=(F_{n_3}^*\otimes I_n)diag(\hat{\mathcal H})^*(F_{n_3}\otimes I_n),
        \end{align*}
        and since $diag(\hat{\mathcal H})$ is Hermitian, we conclude that $\mathtt{bcirc_z}(\mathcal H)$ is Hermitian. Then, by Theorem \ref{newrelation}, Theorem \ref{Hermitian} and equation \eqref{zhao2.1}, we obtain  $ \mathcal A=\mathcal U*_Q\mathcal H$, where $\mathcal H$ is Hermitian and f-positive semidefinite. 
        
        On the other hand, by equation \eqref{diagpsd} and Theorem \ref{zdc}, we obtain 
        \begin{align*}
            x^*(F_{n_3}\otimes I_n)\mathtt{bcirc_z}(\mathcal H)(F_{n_3}^*\otimes I_n)x\geq0.
        \end{align*}
        Let $(F_{n_3}^*\otimes I_n)x=y$, and thus $y^*=x^*(F_{n_3}\otimes I_n)$. 
Since $x$ is arbitrary, $y$ is arbitrary as well. Consequently,
\[
y^{*}\,\mathtt{bcirc_z}(\mathcal H)\,y \ge 0,
\]
where $\mathtt{bcirc_z}(\mathcal H)$ is Hermitian and positive semidefinite. 
    \end{proof}
\end{theorem}

Based on the preceding theorem, we next present an algorithm for computing the QT-polar decomposition, stated in Algorithm \ref{QTpolar}.
\begin{algorithm}[H] \caption{QT-polar decomposition}\label{QTpolar}
\begin{algorithmic}[1]
\smallskip
\Input Third-order quaternion tensor $\mathcal{A} \in \mathbb{Q}^{n \times n \times n_3}$
\Output Third-order unitary tensor $\mathcal{U} \in \mathbb{Q}^{n \times n \times n_3}$, and Hermitian and f-positive semidefinite tensor  $\mathcal{H} \in \mathbb{Q}^{n \times n \times n_3}$.
\State $\hat{\mathcal A}=fftq(\mathcal A,[\ ],3)$
\For {$i=1,2,...,n_3$ }
 $[\hat{\mathcal{U}}(:,:,i),\hat{\mathcal{H}}(:,:,i)]=Polar(\hat{\mathcal A}(:,:,i))$
\EndFor
\State 
$\mathcal{U}=ifftq(\hat{\mathcal{U}},[\ ],3), \quad \mathcal{H}=ifftq(\hat{\mathcal{H}},[\ ],3) \quad$
\end{algorithmic}
\end{algorithm}

To verify the accuracy of Algorithm~\ref{QTpolar} for computing the QT-polar decomposition, we present the following example.
\begin{example}\label{expPolar}
Let $\mathcal{A} = \mathcal{A}_{\mathbf{d}} + \mathbf{j}\mathcal{A}_{\mathbf{c}}  \in \mathbb{Q}^{3 \times 3 \times2}$ with entries
\[\mathcal{A}_{\mathbf{d}}(:,:,1) = \begin{bmatrix} \ 8 & 3+8i & 6+6i\ \\ \ 2 & 2+10i & 5+5i\ \\ \ 10+5i & 2+1i & 3\ \end{bmatrix}, \quad \mathcal{A}_{\mathbf{d}}(:,:,2) = \begin{bmatrix} \ 9+3i & 10+3i & 8+6i\ \\ \ 6+1i & 3+5i & 4+2i\ \\ \ 6+8i & 8+1i & 6+7i\ \end{bmatrix},\] 
and
\[\mathcal{A}_{\mathbf{c}}(:,:,1) = \begin{bmatrix} \ 7 & 0-8i & 1-2i\ \\ \ 8-4i & 2-4i & 9-2i\ \\ \ 4-2i & 10-10i & 5-1i\ \end{bmatrix}, \quad \mathcal{A}_{\mathbf{c}}(:,:,2) = \begin{bmatrix} \ 10-1i & 1-6i & 8-6i\ \\ \ 0-9i & 10-1i & 8-3i\ \\ \ 4-6i & 0-9i & 9-5i\ \end{bmatrix}.\]
Then, applying Algorithm \ref{QTpolar}, we obtain a unitary tensor $\mathcal{U}=\mathcal{U}_{\mathbf{d}} + \mathbf{j}\mathcal{U}_{\mathbf{c}}$, where
\[\mathcal{U}_{\mathbf{d}}(:,:,1) = \begin{bmatrix} \ 0.0977-0.0957i & -0.0911+0.2117i & 0.0162+0.2244i\ \\ \ -0.0457-0.0162i & 0.002+0.3894i & 0.2925+0.0834i\ \\ \ 0.1745+0.214i & 0.0278-0.1106i & -0.0672-0.2074i\ \end{bmatrix},\]
     \[\mathcal{U}_{\mathbf{d}}(:,:,2) = \begin{bmatrix} \ 0.1729+0.0002i & 0.3501-0.1514i & 0.1663+0.1022i\ \\ \ 0.1367+0.0786i & -0.0408+0.1717i & 0.206-0.1925i\ \\ \ 0.1029+0.2674i & 0.1983-0.1068i & 0.033+0.3279i\ \end{bmatrix},\]
     and
\[\mathcal{U}_{\mathbf{c}}(:,:,1) = \begin{bmatrix} \ 0.122+0.0194i & -0.1388-0.182i & -0.1728-0.0508i\ \\ \ 0.2798+0.0329i & -0.0719+0.0795i & 0.1302+0.0373i\ \\ \ 0.1548+0.022i & 0.427-0.2406i & -0.0389+0.2302i\ \end{bmatrix},\]
     \[\mathcal{U}_{\mathbf{c}}(:,:,2) = \begin{bmatrix} \ 0.4842+0.1635i & 0.0403+0.1116i & 0.229-0.4861i\ \\ \ -0.3868-0.4185i & 0.2866+0.1137i & 0.2683-0.1304i\ \\ \ -0.1303-0.2095i & -0.1532-0.3437i & 0.3163-0.0106i\ \end{bmatrix},\]
     and a Hermitian and f-positive semidefinite tensor $\mathcal{H} = \mathcal{H}_{\mathbf{d}} + \mathbf{j}\mathcal{H}_{\mathbf{c}},$ where
\[\mathcal{H}_{\mathbf{d}}(:,:,1) = \begin{bmatrix} \ 21.4335 & 3.0005-1.9322i & 9.0186+1.9319i\ \\ \ 3.0005+1.9322i & 23.4129 & 9.018+2.7682i\ \\ \ 9.0186-1.9319i & 9.018-2.7682i & 18.784\ \end{bmatrix},\]
     \[\mathcal{H}_{\mathbf{d}}(:,:,2) = \begin{bmatrix} \ 10.3414 & 7.0911+0.8509i & 6.3135+3.4139i\ \\ \ 7.0911-0.8509i & 7.0046 & 9.124+1.3119i\ \\ \ 6.3135-3.4139i & 9.124-1.3119i & 6.554\ \end{bmatrix},\]
     and
\[\mathcal{H}_{\mathbf{c}}(:,:,1) = \begin{bmatrix} \ 0 & 1.8096-3.7384i & 0.3446-0.4488i\ \\ \ -1.8096+3.7384i & 0 & 1.648+3.1267i\ \\ \ -0.3446+0.4488i & -1.648-3.1267i & 0\ \end{bmatrix},\]
     \[\mathcal{H}_{\mathbf{c}}(:,:,2) = \begin{bmatrix} \ 0 & -1.2903-3.9939i & 0.397-0.9932i\ \\ \ 1.2903+3.9939i & 0 & -1.4261+2.9459i\ \\ \ -0.397+0.9932i & 1.4261-2.9459i & 0\ \end{bmatrix},\]
     which satisfy $\mathcal{A}=\mathcal{U}*_Q\mathcal{H}$.
\end{example}

\begin{remark}
    Similarly, we can demonstrate the left polar decomposition, that is
    $\mathcal A=\mathcal K*_Q\mathcal W,$ where $\mathcal W$ is a unitary quaternion tensor and $\mathcal K$ is a Hermitian and f-positive semidefinite quaternion tensor. Unless otherwise indicated, QT-polar decomposition herein refers to the QT-right polar decomposition.  
\end{remark}

\section{Large-scale experiments and Scalability analysis}

This chapter presents some large-scale tests and scalability analysis. All computations were performed in MATLAB R2022b using the QTFM quaternion toolbox and the Tensorlab package, and executed on a machine equipped with a 13th Gen Intel(R) Core(TM) i5-13500HX (2.50 GHz) processor and 16 GB of RAM.

We first propose the $z$-block circulant structure-accelerated inversion algorithm ($\mathtt{bcirc_z}$-inv) to solve the problem of computing the inverse of a $z$-block circulant matrix.
For comparison, the built-in quaternion matrix inverse function \texttt{inv} from the
\textit{MATLAB} toolbox \textit{QTFM} is employed.

\begin{algorithm}[H] \caption{$\mathtt{bcirc_z}$-inv}\label{QTinverse}
\begin{algorithmic}[1]
\smallskip
\Input The $z$-block circulant matrix $A=\mathtt{bcirc_z}{(\mathcal A)}$
\Output $A^{-1}$ is the inverse of $\mathtt{bcirc_z}{(\mathcal A)}$
\State Apply the operator $\mathtt{bcirc_z}^{-1}$ to $A$ to obtain a third-order quaternion tensor
    \[
    \mathcal A = \mathtt{bcirc_z}^{-1}(A)
    \in \mathbb{Q}^{n_1 \times n_1 \times n_3}.
    \]

\State Perform the quaternion tensor Fourier transform along the third mode:
    \[
    \widehat{\mathcal A}
    = \mathtt{fftq}(\mathcal A, [~], 3).
    \]
\State  Compute the inverse of each frontal slice of $\widehat{\mathcal A}$, and denote the resulting tensor by
    $\widehat{\mathcal A}^{-1}$.
\State Apply the inverse quaternion tensor Fourier transform:
    \[
    \mathcal A^{-1}
    = \mathtt{ifftq}(\widehat{\mathcal A}^{-1}).
    \]
\State Finally, apply the $\mathtt{bcirc_z}$ operator to $\mathcal A^{-1}$ to obtain the inverse matrix
    \[
    A^{-1} = \mathtt{bcirc_z}(\mathcal A^{-1}).
    \]
\end{algorithmic}
\end{algorithm}

To verify the effectiveness of the proposed method, we generate a quaternion
$z$-block circulant matrix and compute its inverse. The elements in the first column of the quaternion $z$-block circulant matrix are randomly generated
from a uniform distribution using the \texttt{randq} function provided by the
\textit{MATLAB} toolbox \textit{QTFM}.
For each matrix dimension, the computational time is recorded and reported in Table~\ref{tab:inverse_comparison}.

\begin{example}\label{bcirczinverse}Let
$A \in \mathbb{Q}^{n_1 n_3 \times n_1 n_3}$
be a $z$-block circulant matrix of quaternion tensor $\mathcal A$. The inverse of $A$ is computed through the Algorithm~\ref{QTinverse}. The accuracy of the computed inverse is assessed by measuring the deviations of both $AB$ and $BA$ from the identity matrix, which is defined as
\begin{equation*}
\mathrm{err}
=
\max\left\{
\|AB - I\|,
\;
\|BA - I\|
\right\}.
\end{equation*}
Here, $\|\cdot\|$ denotes the Frobenius norm of a quaternion matrix, and the matrix $B$ denotes the inverse matrix computed by the corresponding computational method.
Recall that the proposed fast method is referred to as the $\mathtt{bcirc_z}$-inv, while the method based on the built-in inverse function is denoted by
$\texttt{inv}$ in Table~\ref{tab:inverse_comparison}.

       \begin{table}[htbp]
    \centering
    \caption{Performance Comparison of $A^{-1}$
Computation for Different Matrix Sizes}
    \label{tab:inverse_comparison}
    \begin{tabular}{c c c c c}
        \toprule
        size  & \multicolumn{2}{c}{time (s)} & \multicolumn{2}{c}{err} \\
      \cmidrule(lr){1-5}
    $(n_1n_3 \times n_1n_3)$    & $\mathtt{bcirc_z}$-inv & \texttt{inv} &$ \mathtt{bcirc_z}$-inv & \texttt{inv} \\
        \midrule
        $9\times9 $& $0.0096$ & $0.0015$ & $1.32e-15$ & $2.22e-14$ \\
       $ 15\times15 $& $0.0043$ & $0.0025$ & $8.61e-15$ & $8.37e-15$ \\
       $ 45\times45$ & $0.0091$ & $0.0076$ & $7.23e-15$ & $2.14e-13$ \\
       $ 75\times75 $& $0.0126$ & $0.0163$ & $9.31e-14$ & $7.55e-13$ \\
       $ 225\times225 $& $0.0456$ & $0.0460$ & $6.19e-13$ & $2.73e-11$ \\
       $ 1125\times1125$ & $0.2062$ & $0.5183$ & $1.38e-12$ & $1.56e-09$ \\
       $ 1875\times1875$ & $0.5017$ & $1.4846$ & $1.65e-11$ & $1.47e-08$ \\
       $ 2500\times2500 $& $0.7693$ & $2.9866$ & $4.36e-11$ & $1.04e-07$ \\
        $3750\times3750 $& $1.5859$ & $8.6190$ & $3.05e-11$ & $7.84e-08$ \\
       $ 5000\times5000 $& $3.0237$ & $18.9334$ & $2.24e-10$ & $3.87e-07$ \\
        \bottomrule
    \end{tabular}
\end{table}

The results presented in Table~\ref{tab:inverse_comparison} clearly demonstrate the computational efficiency and numerical accuracy of the proposed 
$\mathtt{bcirc_z}$-inv method compared with the conventional quaternion matrix inverse function \texttt{inv} provided by the QTFM toolbox.

In terms of computational efficiency, the $\mathtt{bcirc_z}$-inv method exhibits
significantly lower time costs as the matrix size increases.
For smaller matrix sizes (e.g., $9\times 9$), the built-in \texttt{inv} function is slightly
faster; however, once the matrix size exceeds $75\times 75$, the
$\mathtt{bcirc_z}$-inv method consistently outperforms it.
Moreover, the performance gap widens substantially at larger scales.
For example, for a $5000\times 5000$ matrix, the $\mathtt{bcirc_z}$-inv method requires
only $3.02$ seconds, whereas the \texttt{inv} function takes $18.93$ seconds, corresponding
to an $84\%$ reduction in runtime.
These results clearly demonstrate the scalability of the proposed Fourier-based approach,
which exploits the block-circulant structure to significantly reduce the computational
complexity.

Regarding numerical accuracy, the $\mathtt{bcirc_z}$-inv method also demonstrates
superior numerical stability.
Across all tested matrix sizes, the error values—measured as the deviation from the
identity matrix—are consistently lower than those obtained using the \texttt{inv} method,
often by several orders of magnitude.
For instance, for a $2500\times 2500$ matrix, the error of the $\mathtt{bcirc_z}$-inv
method is $4.36\times 10^{-11}$, whereas the \texttt{inv} method yields an error of
$1.04\times 10^{-7}$.
These results indicate that the proposed method maintains higher numerical precision,
particularly for large-scale problems, where numerical errors are more likely to
accumulate.

In summary, the $\mathtt{bcirc_z}$-inv method provides an efficient, numerically stable,
and scalable framework for computing the inverse of $\mathtt{bcirc_z}(\cdot)$.
These advantages make the proposed approach particularly well-suited for high-dimensional
quaternion-valued computations, with potential applications in areas such as signal
processing, computer graphics, and robotics.
\end{example}

To assess the computational performance and scalability of the proposed QT-polar decomposition (Algorithm~\ref{QTpolar}), we conduct numerical experiments on quaternion tensors of different dimensions. The corresponding results are reported in Table
~\ref{tab:qt_polar_performance}.
\begin{example}
Let $\mathcal A\in \mathbb{Q}^{n\times n\times n_3}$. As the dimension of the tensor $\mathcal A$ increases, the computation time increases steadily, ranging from $0.0262$ seconds for the smallest $5 \times 5 \times 5$ tensor to $373.0626$ seconds for the largest $300 \times 300 \times 100$ tensor. More significantly, the residual error remains highly stable across all tensor dimensions, maintaining an exceptionally low level from $2.3631 \times 10^{-14}$ for the $5 \times 5 \times 5$ tensor to $8.6245 \times 10^{-10}$ for the $300 \times 300 \times 100$ tensor. 

This trend indicates that the algorithm preserves high numerical precision even for large-scale tensors, with the residual error increasing only gradually by approximately four orders of magnitude as the tensor size grows by several orders of magnitude. The consistent maintenance of low residual errors across all tested configurations highlights the strong numerical stability of the proposed QT-polar decomposition, indicating its suitability for applications requiring high-precision quaternion tensor computations.
\begin{table}[htbp]
\centering
\small
\caption{Performance of the proposed QT-polar decomposition for different tensor sizes}
\setlength{\tabcolsep}{45pt}
\label{tab:qt_polar_performance}
\begin{tabular}{ccc}
\toprule
Size of $\mathcal A$ & \multicolumn{1}{c}{Time (s)} & \multicolumn{1}{c}{Residual Error} \\
\midrule
$5 \times 5 \times 5$ & 0.0262 & $2.3631 \times 10^{-14}$ \\
$5 \times 5 \times 20$ & 0.1199 & $4.9914 \times 10^{-14}$ \\
$5 \times 5 \times 50$ & 0.2516 & $8.6008 \times 10^{-14}$ \\
$5 \times 5 \times 100$ & 0.5007 & $1.2792 \times 10^{-13}$ \\
$20 \times 20 \times 5$ & 0.0956 & $3.4114 \times 10^{-13}$ \\
$20 \times 20 \times 20$ & 0.3809 & $8.5482 \times 10^{-13}$ \\
$20 \times 20 \times 50$ & 0.9557 & $1.3154 \times 10^{-12}$ \\
$20 \times 20 \times 100$ & 1.8294 & $1.7018 \times 10^{-12}$ \\
$50 \times 50 \times 5$ & 0.3274 & $2.1320 \times 10^{-12}$ \\
$50 \times 50 \times 20$ & 1.1800 & $6.3505 \times 10^{-12}$ \\
$50 \times 50 \times 50$ & 2.8643 & $8.3917 \times 10^{-12}$ \\
$50 \times 50 \times 100$ & 6.0996 & $1.2741 \times 10^{-11}$ \\
$150 \times 150 \times 5$ & 2.8986 & $3.6869 \times 10^{-11}$ \\
$150 \times 150 \times 20$ & 10.5818 & $6.6374 \times 10^{-11}$ \\
$150 \times 150 \times 50$ & 26.8217 & $1.1123 \times 10^{-10}$ \\
$150 \times 150 \times 100$ & 53.2343 & $1.6079 \times 10^{-10}$ \\
$300 \times 300 \times 5$ & 15.6811 & $1.4235 \times 10^{-10}$ \\
$300 \times 300 \times 20$ & 61.4296 & $3.0754 \times 10^{-10}$ \\
$300 \times 300 \times 50$ & 169.2669 & $4.9957 \times 10^{-10}$ \\
$300 \times 300 \times 100$ & 373.0626 & $8.6245 \times 10^{-10}$ \\
\bottomrule
\end{tabular}
\end{table}

Finally, a scalability analysis is conducted, with the corresponding results shown in Figure~\ref{fig:Polar_analysis}. The figure illustrates how the computational efficiency and numerical error of Algorithm~\ref{QTpolar} vary with the tensor size.

\begin{figure}[htbp]
    \centering
    \includegraphics[width=0.45\textwidth]{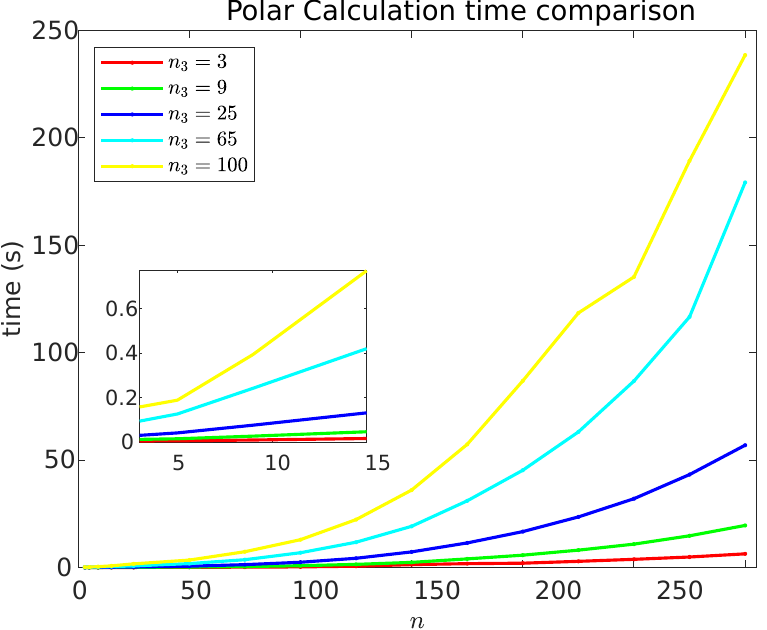}
    \qquad
    \includegraphics[width=0.433\textwidth]{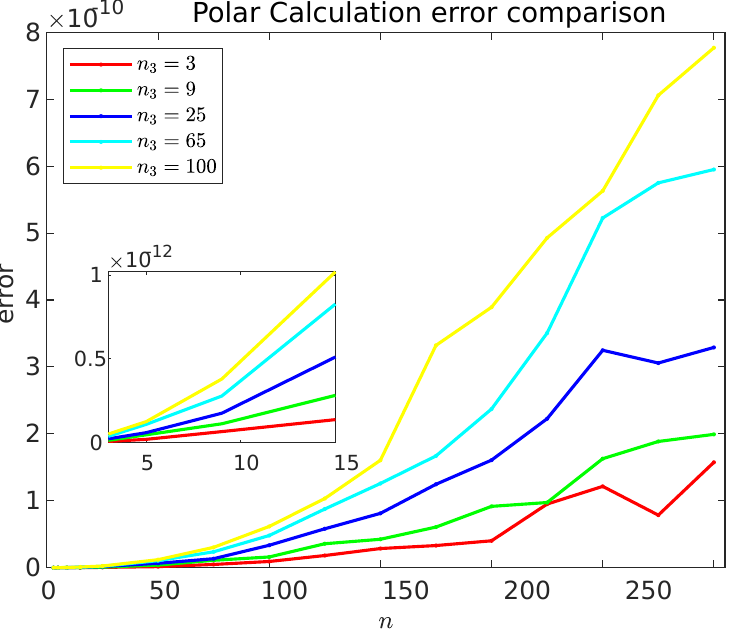}
    \caption{Scalability analysis of QT-polar decomposition algorithm}
    \label{fig:Polar_analysis}
\end{figure}
\end{example}

\section{Applications in Tikhonov-Regularized Model and Video Rotation}

In this section, large-scale testing and comparison of the Tikhonov-regularized model are performed using Theorem \ref{AB} and Algorithm \ref{QTinverse}. Subsequently, a practical example is provided to demonstrate video rotation via the QT-polar decomposition.

Tikhonov regularization is a basic technique for solving ill-posed inverse problems by introducing a penalty term to stabilize the solution. First, we describe the mathematical formulation of Tikhonov regularization in the quaternion domain, which is particularly useful for multidimensional signal processing applications such as color image restoration and vector field analysis.
 
\begin{example}
Consider a linear inverse problem in which a quaternion matrix represents the forward operator. Let $\mathbf{B} \in \mathbb{Q}^{n \times n}$ denote the system matrix, $\mathbf{x} \in \mathbb{Q}^n$ the unknown signal vector, and $\mathbf{b} \in \mathbb{Q}^n$ the observed data vector. The forward model is given by
\begin{equation*}
\mathbf{b} = \mathbf{B}\mathbf{x} + \boldsymbol{\eta},
\end{equation*}
where $\boldsymbol{\eta} \in \mathbb{Q}^n$ represents additive noise.

To address the ill-posedness of the inverse problem, we employ the classical Tikhonov regularization approach. The regularized solution is obtained by minimizing the following objective function:
\begin{equation}
\min_{\mathbf{x}} \left\{ \|\mathbf{B}\mathbf{x} - \mathbf{b}\|_2^2 + \lambda^2 \|\mathbf{x}\|_2^2 \right\},
\label{mini}
\end{equation}
where $\lambda > 0$ is the regularization parameter that balances the trade-off between data fidelity (the first term) and solution regularization (the second term).

The minimization problem in \eqref{mini} admits a closed-form solution by setting the gradient of the objective function to zero. Expanding the objective function yields
\begin{equation*}
J(\mathbf{x}) = (\mathbf{B}\mathbf{x} - \mathbf{b})^H(\mathbf{B}\mathbf{x} - \mathbf{b}) + \lambda^2 \mathbf{x}^H\mathbf{x}.
\end{equation*}
Taking the gradient with respect to $\mathbf{x}$ and setting it to zero yields
\begin{equation*}
\nabla_{\mathbf{x}} J(\mathbf{x}) = 2\mathbf{B}^H(\mathbf{B}\mathbf{x} - \mathbf{b}) + 2\lambda^2 \mathbf{x} = \mathbf{0}.
\end{equation*}
Rearranging terms leads to the associated normal equation
\begin{equation*}
\left( \mathbf{B}^H \mathbf{B} + \lambda^2 \mathbf{I} \right) \mathbf{x} = \mathbf{B}^H \mathbf{b}.
\end{equation*}
Consequently, the solution can be expressed explicitly as
\begin{equation*}
\mathbf{x}_{\text{direct}} = ( \mathbf{B}^H \mathbf{B} + \lambda^2 \mathbf{I})^{-1}\mathbf{B}^H \mathbf{b}.
\end{equation*}

Although the direct matrix inversion method provides an exact solution to the Tikhonov regularized system, its computational complexity of $\mathcal{O}(n^3)$ becomes prohibitive for large-scale testing, therefore practical only for small- to moderate-scale problems.

The motivation is to introduce the system matrix $\mathbf{B}$ that possesses a $z$-block circulant structure $\mathtt{bcirc}_z(\mathcal B)$ and exploit this special algebraic structure to accelerate the computation of the matrix inverse, thereby making large-scale experimental evaluations of the Tikhonov-regularized system using direct methods feasible. 
Given that $$\mathbf{B}^H \mathbf{B} + \lambda^2 \mathbf{I}=\text{$\mathtt{bcirc_z}$}(\mathcal{B}^H*_Q \mathcal{B} + \lambda^2 \mathcal{I}),$$ which is obtained by Theorem \ref{AB} and inherits the $z$-block circulant structure. Then we can compute its inverse efficiently using Example \ref{bcirczinverse}.

The computational efficiency and numerical accuracy of the algorithm $\mathtt{bcirc_z}$-inv were evaluated through comprehensive numerical experiments. The algorithm was implemented in MATLAB and tested on a range of problem sizes, with its performance compared against \texttt{inv}.

Experiments were conducted with various tensor dimensions, where the parameter $m$ represents the block size and $q$ denotes the number of channels. For simplicity, let  $\mathbf{A} = \mathbf{B}^H \mathbf{B} + \lambda^2 \mathbf{I}\in \mathbb{Q}^{n\times n}$, where $n=mq$. And $\lambda = 0.5$ was used as the regularization parameter. The execution times and errors for both $\mathtt{bcirc_z}$-inv and \texttt{inv} are presented in Fig.~\ref{fig:exp_bcirc_z_inv}.

For small-scale testing, both methods perform comparably, with \texttt{inv} occasionally running slightly faster due to its lower overhead. However, as the testing size increases, $\mathtt{bcirc_z}$-inv demonstrates significantly superior scaling behavior.
\begin{figure}[htbp]
    \centering
    \includegraphics[width=0.45\textwidth]{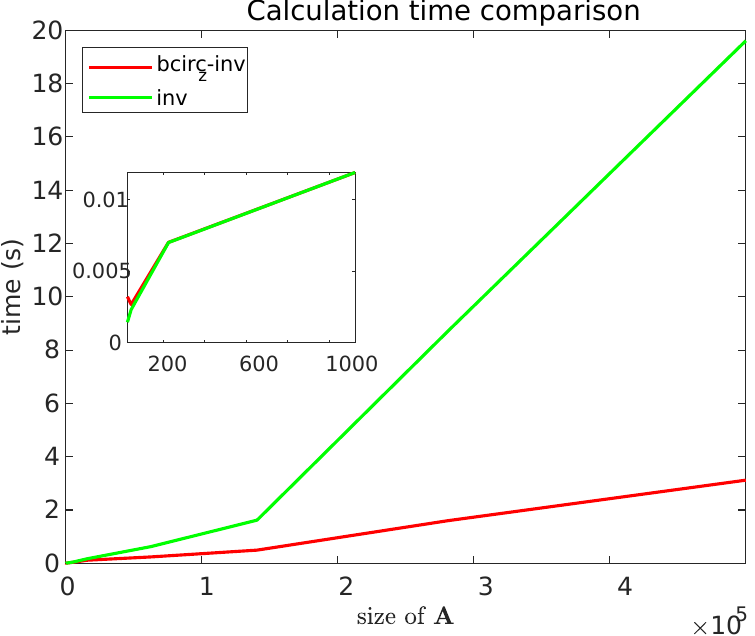}
    \qquad
    \includegraphics[width=0.45\textwidth]{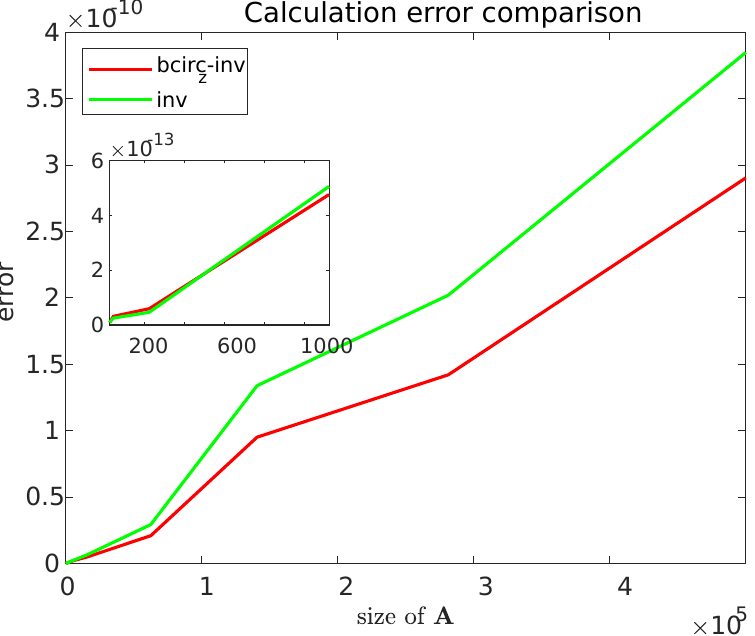}
    \caption{Scalability analysis of computing $\mathbf A^{-1}$ with two methods}
    \label{fig:exp_bcirc_z_inv}
\end{figure}
At moderate problem sizes $(n = 225)$, $\mathtt{bcirc_z}$-inv begins to show measurable advantages, achieving execution times approximately 20\% faster than those of \texttt{inv}. This performance gap widens substantially as the problem dimension increases. Despite the computational acceleration, $\mathtt{bcirc_z}$-inv maintains comparable numerical accuracy. The residual norms of both methods are on the order of $10^{-11}$ to $10^{-14}$ for smaller problems and $10^{-10}$ to $10^{-11}$ for larger ones. In the largest case $(n = 5000)$, $\mathtt{bcirc_z}$-inv yields a residual of $2.90 \times 10^{-10}$, which is still slightly lower than the $3.85 \times 10^{-10}$ obtained by \texttt{inv}. The performance gain stems from exploiting the $z$-block circulant structure, which reduces the computational complexity from $\mathcal{O}(n^3)$ $(\mathtt{inv})$ to $\mathcal{O}(m^3 q + m^2 q^2)$ $(\mathtt{bcirc_z}-inv)$, where $n = m \times q$. 

In summary,  $\mathtt{bcirc_z}$-inv exploits the special algebraic structure of the system matrix to achieve computational efficiency while maintaining numerical accuracy. This approach reduces the overall computational cost and enables the efficient solution of large-scale testing. Therefore, breaking the $\mathtt{inv}$ limitation to small- to moderate-scale testing.  For problems with a $z$-block circulant structure and large dimensions, $\mathtt{bcirc_z}$-inv provides a practical and efficient solution to the Tikhonov-regularized inverse problem, making it suitable for applications in color image processing and multi-channel signal analysis.
\end{example}

Polar decomposition has been widely applied in continuum mechanics and quantum mechanics, specifically in computer graphics, such as image 
rotation~\cite{KT1992}, face recognition~\cite{MKIET2011}, and image watermarking~\cite{QSSVC405}. 
\begin{example}
In essence, the polar decomposition factors a shear-scaling matrix into the product of a rotation 
matrix and a scaling matrix. For matrices over the complex field, this corresponds to expressing a 
matrix as the product of a unitary matrix and a positive definite Hermitian matrix. In this 
framework, the unitary (or rotation) matrix $Q$ captures the global rotational component, whereas 
the positive definite Hermitian matrix $H$ characterizes the anisotropic scaling along the 
principal directions determined by its eigenvectors.

Due to its decomposition property that the polar decomposition separates deformation into pure rotation and pure scaling, we extend its application to quaternion tensors to handle the problems of video rotation, which allows for rotating different RGB channels by different angles at different frames of a video.

It is worth noting that video rotation is not only a correction tool but also a means of creative expression. In science fiction or horror films, rotating different color channels to varying degrees can depict a character's sense of loss in an alternate dimension or gravity. Traditional glitch art often involves misalignment of color channels. However, upgrading this misalignment from simple translation to rotation makes the glitch effect more dynamic and complex, no longer a cold digital error but with an organic, vortex-like texture. This technique is widely used in cyberpunk-themed films, music videos (especially electronic music, psychedelic rock), and artistic short films depicting dreams.

In what follows, we introduce an application example of using polar decomposition for video rotation. Given a rotation-scaling tensor $ \mathcal C \in \mathbb{Q}^{2\times 2\times p} $, we use QT-polar decomposition to obtain $\mathcal  C =\mathcal  U*_Q\mathcal  H $, where $ \mathcal U $ is a unitary tensor, and $ H $ is a Hermitian and f-positive semidefinite tensor. The resulting $ \mathcal U $ is a tensor containing rotation information for the R, G, and B channels of the video

We apply $ \hat {\mathcal U}^{(n)}= \text{qfft}(\mathcal U^{(n)}, [\ ], 3) $, for $ n = 1, 2, \dots, p $,
\[
\hat {\mathcal U}^{(n)} = 
\begin{pmatrix}
\cos\alpha_n & 0 \\
0 & \cos\gamma_n
\end{pmatrix} i + 
\begin{pmatrix}
\sin\alpha_n \cos\beta_n & 0 \\
0 & \sin\gamma_n \cos\beta_n
\end{pmatrix} j + 
\begin{pmatrix}
\sin\alpha_n \sin\beta_n & 0 \\
0 & \sin\gamma_n \sin\beta_n
\end{pmatrix} k
\]

Thus, by obtaining the rotation parameters $\alpha, \beta, \gamma$, we conducted experiments in which the R, G, and B channels of different video frames were rotated using both identical and distinct angles, thereby producing a rotated video. The specific results are presented in Table~\ref{tab_spin}.

\begin{table}[htbp]
\centering
\begin{threeparttable}
\caption{Comparison of Image Rotation Results}
\label{tab_spin}
\begin{tabular}{ccccccc}
\toprule
frame &\multicolumn{2}{c}{120} & \multicolumn{2}{c}{230} & \multicolumn{2}{c}{340}  \\
 Rotation method &image & Angle  &image & Angle &image & Angle \\
\midrule
\multirow{1}{*}[16pt]{Original}& 
\includegraphics[width=0.15\textwidth]{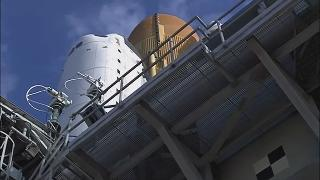} & \multirow{1}{*}[23pt]{\shortstack{0\\0\\0}}&
\includegraphics[width=0.15\textwidth]{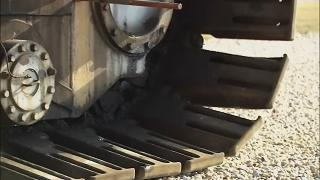} & \multirow{1}{*}[23pt]{\shortstack{0\\0\\0}}&
\includegraphics[width=0.15\textwidth]{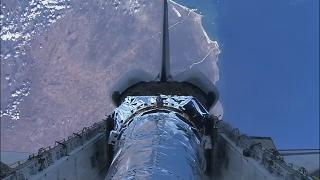} & \multirow{1}{*}[23pt]{\shortstack{0\\0\\0}}\\
\multirow{1}{*}[16pt]{Same linearity\tnote{a}}& 
\includegraphics[width=0.15\textwidth]{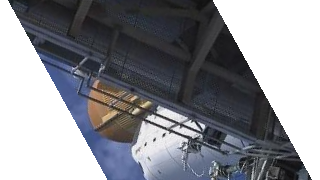} & \multirow{1}{*}[23pt]{\shortstack{119\\119\\119}}&
\includegraphics[width=0.15\textwidth]{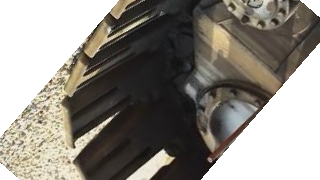} & \multirow{1}{*}[23pt]{\shortstack{220\\229\\229}}&
\includegraphics[width=0.15\textwidth]{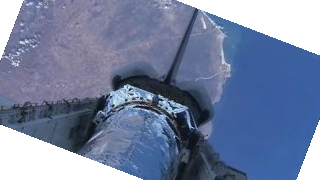} & \multirow{1}{*}[23pt]{\shortstack{339\\339\\339}}\\
\multirow{1}{*}[16pt]{Different linearity\tnote{b}}& 
\includegraphics[width=0.15\textwidth]{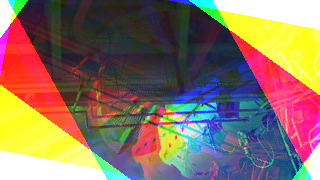} & \multirow{1}{*}[23pt]{\shortstack{119\\144\\169}}&
\includegraphics[width=0.15\textwidth]{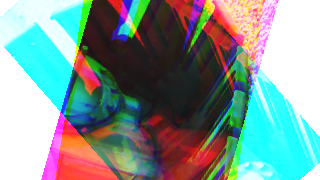} & \multirow{1}{*}[23pt]{\shortstack{131\\74\\81}}&
\includegraphics[width=0.15\textwidth]{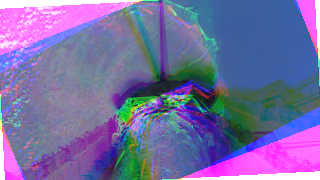} & \multirow{1}{*}[23pt]{\shortstack{21\\3\\28}}\\
\multirow{1}{*}[16pt]{Fixed step\tnote{c}}& 
\includegraphics[width=0.15\textwidth]{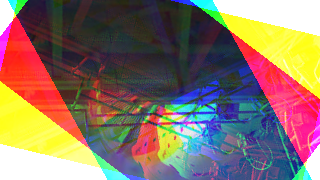} & \multirow{1}{*}[23pt]{\shortstack{119\\142.8\\166.6}}&
\includegraphics[width=0.15\textwidth]{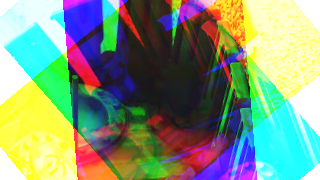} & \multirow{1}{*}[23pt]{\shortstack{131\\94.8\\39.4}}&
\includegraphics[width=0.15\textwidth]{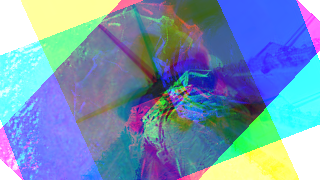} & \multirow{1}{*}[23pt]{\shortstack{21\\46.8\\114.6}}\\
\multirow{1}{*}[16pt]{Sine wave phase\tnote{d}}& 
\includegraphics[width=0.15\textwidth]{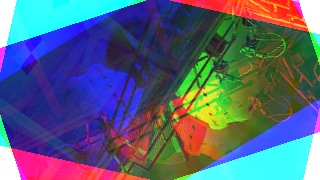} & \multirow{1}{*}[23pt]{\shortstack{167.2\\198.1\\30.9}}&
\includegraphics[width=0.15\textwidth]{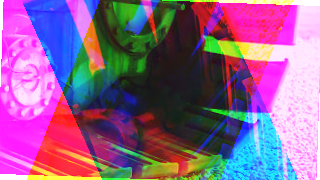} & \multirow{1}{*}[23pt]{\shortstack{114.2\\358.0\\68.0}}&
\includegraphics[width=0.15\textwidth]{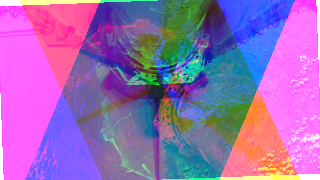} & \multirow{1}{*}[23pt]{\shortstack{67.2\\182.1\\115.0}}\\
\multirow{1}{*}[16pt]{Logarithmic growth\tnote{e}}& 
\includegraphics[width=0.15\textwidth]{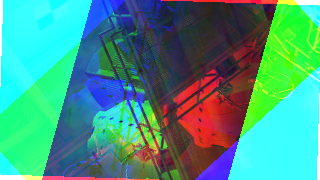} & \multirow{1}{*}[23pt]{\shortstack{178.2\\74.3\\42.2}}&
\includegraphics[width=0.15\textwidth]{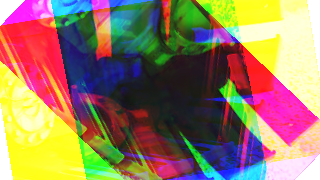} & \multirow{1}{*}[23pt]{\shortstack{101.3\\316.1\\9.3}}&
\includegraphics[width=0.15\textwidth]{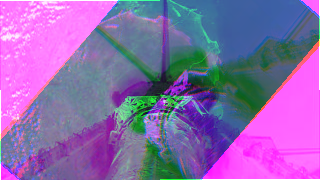} & \multirow{1}{*}[23pt]{\shortstack{50.0\\0.2\\48.1}}\\
\bottomrule
\end{tabular}
\begin{tablenotes}
\item[a] Each channel rotates through a fixed range with constant angular velocity. And each color channel has the same initial angle.
\item[b] Each channel rotates through a fixed range with constant angular velocity. And each color channel has a different initial angle.
\item[c] Each channel rotates with a constant step per frame, leading to different angular velocities. And each color channel has the same initial angle.
\item[d] Each channel rotates with a sinusoidal pattern with a fixed phase shift between channels. And each color channel has the same initial angle.
\item[e] Each channel rotates with a logarithmic function, creating a non-linear change over time. And each color channel has the same initial angle.
\end{tablenotes}
\end{threeparttable}
\end{table}

For quantitative evaluation, we adopt the following metrics. Specifically, we compute the mean Temporal Consistency (TC${\text{mean}}$), the mean Color Consistency (CC${\text{mean}}$), and the standard deviation of Temporal Consistency (TC$_{\text{std}}$), which are defined as follows.
\begin{align*}
\text{TC}_{\text{mean}} &= \frac{1}{N-1} \sum_{i=1}^{N-1} \left[1 - \frac{\mathbb{E}[|\mathcal G_{i+1} - \mathcal G_i|]}{255}\right], \\
\text{CC}_{\text{mean}} &= \frac{1}{N} \sum_{i=1}^{N} \frac{1}{1 + \bar{\sigma}_i/255}, \\
\text{TC}_{\text{std}} &= \sqrt{\frac{1}{N-2} \sum_{i=1}^{N-1} (\text{TC}_i - \text{TC}_{\text{mean}})^2},
\end{align*}
where $\mathcal{G}_i$ denotes the grayscale image of the $i$-th frame, $\mathbb{E}[\cdot]$ represents the mathematical expectation (i.e., the mean), $\bar{\sigma}_i$ is the average standard deviation of the three color channels of the $i$-th frame, and $N$ denotes the frame rate of the video. The evaluation results are summarized in Table~\ref{tab_assess_spin}.

By testing several different rotation methods, we observe two key aspects: (1) The data in the temporal consistency standard deviation of the metrics also achieves favorably, which demonstrates the stability of the QT-polar decomposition. (2) The data in both temporal consistency mean and color consistency mean perform favorably, which demonstrates that the QT-polar decomposition preserves frame-to-frame coherence while maintaining stability and faithful color reproduction. 
\begin{table}
\centering
\caption{Assessment of Rotation Results}
\label{tab_assess_spin}
\begin{tabular}{cccc}
\toprule
&temporal consistency mean& temporal consistency std&color consistency mean\\
\midrule
Same linearity&0.9652&0.0172&0.7708\\
Different linearity&0.9648&0.0202&0.7701\\
Fixed step size&0.9625&0.0190&0.7713\\
Sine wave phase&0.9440&0.0236&0.7815\\
Logarithmic growth&0.9707&0.0206&0.7741\\
\bottomrule
\end{tabular}
\end{table}

\end{example}

\section*{Funding}
\begin{itemize}
\item The first author was supported by the National Natural Science Foundation of China (NSFC) (No. 11901079), and China Postdoctoral Science Foundation (No. 2021M700751), and the Scientific and Technological Research Program Foundation of Jilin Province (No. JJKH20190690KJ; No. JJKH20220091KJ; No. JJKH20250851KJ).
\item The fourth author is supported by the Ministry of Science, Technological Development and Innovations, Republic of Serbia (No. 451-03-137/2025-03/200124).
\end{itemize}

\section*{Conflict of Interest}
The authors declare that they have no potential conflict of interest.

\section*{Data Availability}
Data sharing is not applicable to this article, as no datasets were generated or analyzed during the current study.


\bibliographystyle{abbrv}

\end{document}